\documentclass[11pt,a4paper]{article}
\usepackage{authblk}
\usepackage{graphicx}
\usepackage{latexsym}
\usepackage{subfigure}
\usepackage{amsmath}
\usepackage{enumerate}
\usepackage{amssymb}
\usepackage[mathscr]{eucal}
\usepackage{amsthm}
\usepackage{color}
\usepackage{array}
\usepackage[small]{caption}
\usepackage{hyperref}
\usepackage[a4paper,text={6.0in,9.0in},centering,
            includefoot,foot=0.6in]{geometry}

\usepackage[normalem]{ulem}

\allowdisplaybreaks[1]

\newtheorem{thm}{Theorem}[section]
\newtheorem{defn}[thm]{Definition}
\numberwithin{equation}{section}

\newcommand{\bT}{\mathbf{T}}
\newcommand{\R}{\mathcal{R}}

\renewcommand{\L}{\mathcal{L}}
\renewcommand{\P}{\mathcal{P}}

\newcommand{\emat}{\end{pmatrix}}

\DeclareMathOperator{\Id}{{\rm Id}}

\begin{document}


\title{A micro/macro parallel-in-time (parareal) algorithm applied to a climate model with discontinuous non-monotone coefficients and oscillatory forcing}

\author[1]{Giovanni Samaey}
\author[2]{Thomas Slawig}
\affil[1]{
\small
Department of Computer Science, K.U. Leuven,
Celestijnenlaan 200A, 3001 Leuven, Belgium, \tt{giovanni.samaey@kuleuven.be}}
\affil[2]{ 
\small Department of Computer Science, Kiel University,  24098 Kiel, Germany, \tt{ts@informatik.uni-kiel.de}
}

\maketitle

\begin{abstract}
We present the application of a micro/macro parareal algorithm for a 1-D energy balance climate model with discontinuous and non-monotone coefficients and forcing terms. The micro/macro parareal method uses a coarse propagator, based on a (macroscopic) 0-D approximation of the underlying (microscopic) 1-D model. We compare the performance of the method using different versions of the macro model, as well as different numerical schemes for the micro propagator, namely an explicit Euler method with constant stepsize and an adaptive library routine. 
We study convergence of the method and the theoretical gain in computational time in a realization on parallel processors. We show that, in this example and for all settings, the micro/macro parareal method converges in fewer iterations than the number of used parareal subintervals, and that a theoretical gain in performance of up to 10 is possible. 
\end{abstract}

\maketitle


\section{Introduction}

Climate simulation is among the most challenging and time-consuming computational tasks, for a number of reasons. First, there is the complexity of the coupled climate system, with interactions between many different components and nonlinearity of many important processes, some of which are not completely understood by now (see e.g. \cite{McgHen14}).
Second, the need for high spatial resolution in global climate models results in a huge dimension of the discretized systems that have to be solved. To reduce the needed computational time, spatial parallelization is a common strategy used in fully coupled high resolution climate models.  Still, for long-time  simulation  runs (e.g., to compute full glacial cycles of  hundreds of thousands of years), spatial parallelization alone is insufficient, since spatial and temporal resolutions are typically coupled via some kind of CFL-type condition. Thus,  such long-time simulation runs are  only feasible with lower complexity models (see, e.g.,  \cite{Ganopolski2017}.  To proceed to long-time \emph{and} high resolution models, an additional parallelization in time becomes very attractive.

In this paper, we present how we generalized and used a micro/macro version of a parallel-in-time algorithm for time integration, called \emph{parareal} \cite{LiMaGa01}, as it was developed in \cite{LeLeSa13}.  
 As an {example} problem for a climate model, we choose an energy balance model (EBM) that describes the evolution of global mean temperature by balancing incoming and outgoing radiation, see, e.g., \cite[Section 3]{McgHen14}). Such models are the simplest way to model the Earth's climate. Incoming radiation is determined by the energy received from the Sun, diminished by a fraction (called albedo) which is reflected by the Earth's atmosphere (e.g., clouds) or surface (e.g., ice).  Outgoing radiation is usually determined by considering the Earth as a perfect black radiating body, for which the Stefan-Boltzmann law gives a relation between the outgoing radiation and the temperature of the body. Since -- due to the greenhouse effect -- the Earth is not perfectly radiating, an emissivity parameter is included in EBMs to take into account that part of the emitted heat is captured in the Earth's atmosphere. 
Here, we regard as the ``original'' or ``micro'' model an EBM in one space dimension, using latitude as the spatial coordinate. 

The micro/macro parareal method will use an approximate ``macro'' model to obtain a numerically fast predictor, which is iteratively corrected by time-parallel simulations using the original 1-D model on different slices of the time interval. The macro model is obtained by {considering only global mean values of temperature. It can be obtained by spatially} averaging {or simplifying} the micro 1-D model. {This} results in a 0-D model that takes the form of an ordinary differential equation (ODE). The resulting macro model is much cheaper to simulation numerically for two reasons. First, due to its low-dimensionality, fewer degrees of freedom need to be accounted for. Second, the macro model only contains the dominant slow time scale, such that larger time steps can be taken.  
We mainly use the 0-D model to accelerate simulations. Nevertheless, 0-D EBMs have their own justification,  for educational purposes (see \cite{McgHen14}), 
since they include most important features of the Earth's energy balance. They model the  Earth as a point in space, and thus all parameters mentioned above 
 (energy received from the Sun as only external forcing of the climate system as well as albedo and emissivity  that enter the differential equations as coefficients)
 represent averaged values over the whole planet. {In  0-D models, nonlinear dependencies of coefficients on the state variable (temperature) usually are restricted to the albedo.}

If only a few parareal iterations are required, the micro/macro parareal algorithm can achieve a significant reduction in the required wall-clock time, compared to a naive simulation using only the micro-simulator by performing the 1-D simulations on different slices of the time interval \emph{in parallel} in each parareal iteration. In that case, one can reach an accuracy that is much higher than that of the macroscopic simulation on the whole time interval of interest, with a wall-clock time that is of the order of the simulation time of a few parareal time slices. In this paper, we study to what extent this potential is realized for climate simulations based on EBMs. In particular, we investigate the ability of the method to deal with oscillatory and abrupt changes in forcing terms, and we study the dependency of the method on the accuracy of the macro model and the number of parareal time slices.

Since its introduction in~\cite{LiMaGa01}, the parareal strategy has been applied to a wide range of problems, including fluid-structure interaction~\cite{farhat2003time}, Navier--Stokes equation simulation~\cite{fischer2005parareal}, and reservoir simulation~\cite{garrido2005}. We refer to ~\cite{maday2002parareal,maday2005parareal} for further analysis, and to \cite{bal2005convergence,staff2005stability} for stability results. In~\cite{gander2007analysis}, the method is reformulated in a more general setting that relates the parareal strategy to earlier time-parallel algorithms, such as multiple shooting (see e.g.~\cite{keller1968numerical,nievergelt1964parallel}) or multigrid waveform relaxation (see e.g.~\cite{lubich1987multi,vandewalle1992efficient}).  The micro-macro parareal method in this paper is a generalization of the method in \cite{LeLeSa13}, in which the micro model was a high-dimensional stiff ODE and the macro model was an approximate, low-dimensional ODE for a limited set of slow degrees of freedom. {In this work, the same principle is followed: we design a parareal method in which the coarse propagator uses a lower-dimensional model than the fine propagator. The specific novelty in this paper is the choice of two two models (a 0-D and a 1-D EBM), from which stems the need to design a specific coupling approach to transfer information between the two levels during the serial step of the parareal iteration. }

Similar micro-macro parareal methods have been considered in the literature. The authors of~\cite{BBK,maday41parareal} consider a singularly perturbed system of ordinary differential equations (ODEs) at the microscopic level and the limiting differential-algebraic equation at the macroscopic level. In~\cite{engblom2009parallel}, a parareal algorithm for multiscale stochastic chemical kinetics is presented, in which the macroscopic level uses the mean-field limiting ODE. In~\cite{mitran2010time}, the parareal algorithm is used with kinetic Monte Carlo at the macroscopic level and molecular dynamics at the microscopic level.

{In the climate community, parallel-in-time methods have been used for simple ODE models,  e.g., for the Lorenz model in \cite{Gander2008}. Classical  1-D energy balance models as the one we studied in our work include diffusion terms that  model  energy transport in space. They do \emph{not} include advection terms. Realistic climate models additionally include mass and momentum balance and are based on the Navier-Stokes equations. Moreover, the corresponding flow problems show strong advection. Whereas parallel-in-time algorithms have good convergence properties for diffusion-dominated equations,  dominant advection  causes problems. Thus,  the application of the micro-macro parareal method to realistic climate models has to be further investigated. Our motivation to study 1-D energy balance models was the presence of highly nonlinear coefficients and, in our setting, a multi-scale forcing in time.}

This text is organized as follows. In Section~\ref{sec:model}, we describe the model problem in its 0-D and 1-D versions. We also describe  how the 1-D coefficients are captured in the 0-D version. In Section \ref{sec:parareal}, we describe the micro/macro parareal algorithm that is the focus of the present paper. In Section \ref{sec:numexp}, we describe in detail the numerical experiments that we performed. We present the numerical results in Section~\ref{sec:results}, discussing both the effect of the choice of coefficients in the 0-D model and the choice of micro and macro time integration methods. We end the paper with a summary and conclusions in Section \ref{sec:discussion}.

\section{Model problem\label{sec:model}}

In this Section, we discuss the energy balance models (EBMs) that will be used throughout the paper.  In the 1-D model, we have a single spatial coordinate $\phi$, which varies from $\phi=0$ at the north pole to $\phi= \pi/2$ at the equator, thus named colatitude. Due to a symmetry assumption (following \cite{Bay91,Ghi75}), only half of the sphere is modeled.   The state variable, temperature, then is a function $T=T(\phi,t)$ of colatitude and time. In this model, spatial redistribution of heat is included via diffusion, with a potentially temperature-dependent diffusion coefficient. Incoming radiation, albedo and thermal capacity depend on the colatitude $\phi$. Albedo and emissivity (which models the greenhouse effect) also depend on temperature, to take into account ice melting or increase of water vapor in the atmosphere, respectively, leading to additional nonlinearities in the model.
By {considering only the global mean temperature}, a 0-D model can be obtained, in which the Earth is seen as a point in space. {Thus also all parameters} (energy received from the Sun, albedo, emissivity) are  averaged over the whole planet {in some sense}. The state variable  then is a scalar function $T$ of time $t$.

For clarity of exposition, we first discuss the 0-D model in Section~\ref{sec:model-0-D}, after which we elaborate on the 1-D model in Section~\ref{sec:model-1-D}.  The choice of the parameters in the 1-D model is discussed in Section~\ref{sec:model-params}, along with the connection between the 0-D and 1-D models.  We conclude this Section with some comments on existence and uniqueness of solutions (Section~\ref{sec:model-exist}).

\subsection{0-D model\label{sec:model-0-D}}
In a 0-D EBM \cite[Section 3.2]{McgHen14}, 
we consider the instantaneous change of the time-dependent temperature $T(t)$ due to the difference between incoming and outgoing radiation.
Any difference between ingoing and~outgoing, energy -- $R_{in}(t)$, resp.,~$R_{out}(t)$ -- induces a temporal change of thermal energy,
\begin{eqnarray*}
\left(4\pi r^{2}h\;c\;\rho\right) \; T'(t) = R_{in}(t)-R_{out}(t)
\end{eqnarray*}
in which  $4\pi r^{2} h$ is the volume of the considered spherical shell, $c$ is the specific heat of the fluid, and $\rho$ the respective density.

The total amount of incoming radiant energy per  time unit for the whole Earth  is given as 
\begin{eqnarray}
\label{eq:Rin}
R_{in}&=&(1-\alpha)\pi r^{2} S
\end{eqnarray}
 with  $\alpha\in[0,1]$ the albedo, $\pi r^{2}$ the area of the Earth's two-dimensional projection  ($r$ being the Earth's radius), and $S\approx 1367\,{\rm Wm}^{-2}$ the amount of energy per second and area (often denoted as solar ``constant'', but it fact not constant due to temporal variations of solar activity). In general, the albedo $\alpha$ may be temperature-dependent. We will specify its value in Section~\ref{sec:model-params}.

The outgoing radiant energy per unit time is given by the Stefan-Boltzmann law,
\begin{eqnarray}
\label{eq:Rout}
R_{out}&=&4\pi r^{2}\;\epsilon\;\sigma T^{4}
\end{eqnarray}
 with $4\pi r^{2}$ the Earth's surface and $\sigma= 5.67\times 10^{-8}{\rm Wm}^{-2}\,{\rm K}^{-4}$ the Stefan-Boltzmann constant, and including the emissivity $\epsilon$, i.e., the fraction of outgoing radiation that is not captured by the atmosphere.  Like the albedo $\alpha$, the emissivity $\epsilon$ usually depends on the temperature, see~Section~\ref{sec:model-params}.  

Since $\pi r^{2}$ cancels out in the balance equation, we obtain  
as resulting ODE 
\begin{eqnarray}
\label{eq:0-D}\label{eq:macro}
T'(t) &=&\frac1C\left((1-\alpha)Q- \epsilon\sigma T(t)^4\right),
\end{eqnarray}
where we introduced the symbol $Q:=S/4$ for notational convenience.

In the easiest case of constant solar insolation $Q:=S/4$, heat capacity $C=hc\rho$, emissivity $\epsilon$, and albedo $\alpha$, a stationary solution can be  computed from \eqref{eq:macro}
as
\begin{eqnarray}
\label{eq:stat0-D}
 T_{stat}&=&\sqrt[4]{\frac{(1-\alpha)Q}{ \epsilon\sigma}}.
\end{eqnarray}

\subsection{1-D model\label{sec:model-1-D}}

The following 1-D model, based on \cite{Bay91,Ghi75} and \cite[Section 10]{GhiChi87}, includes variation of the temperature in the latitudinal direction
 $\phi \in [0,\pi/2]$, i.e., we model a half sphere from north pole to equator. To this end, we add {a} diffusive term of the form $\nabla\cdot (k{(\phi,T)}\nabla T(\phi,t))$. Using spherical coordinates and taking into account that colatitude is the only spatial coordinate, this term reduces to
\begin{eqnarray}
\label{eq:diffterm}
\nabla\cdot (k{(\phi,T)}\nabla T(\phi,t))&=&\frac{1}{\sin\phi}\frac{\partial}{\partial \phi}\left(k{(\phi,T)}\sin \phi\frac{\partial T(\phi,t)}{\partial \phi}\right).
\end{eqnarray}

The coefficient $k$ may, in principle, depend on space and temperature. We also include a dependence on the spatial coordinate for the solar radiation $Q$, the albedo $\alpha$, and the heat capacity $C$ a on the spatial coordinate $\phi$. 
In the most general setting, the complete 1-D equation of energy balance gives:
\begin{eqnarray}
\label{eq:1}
C(\phi,T)\frac{\partial T}{\partial t}=
\frac{1}{\sin\phi}\frac{\partial}{\partial \phi}\left(k(\phi,T)\sin \phi\frac{\partial T}{\partial \phi}\right)+\left[1-\alpha(\phi,T)\right]Q(\phi)
 - \epsilon(T)\sigma T^4.
\end{eqnarray}

We use homogeneous Neumann boundary conditions
\begin{eqnarray*}
\frac{\partial T}{\partial \phi}(\phi,t)&=&0\quad\text{for }\phi\in\left\{0,\frac\pi2\right\},t\ge0,
\end{eqnarray*}
where the condition at $\phi=\pi/2$ (equator)  comes from the  assumed symmetry.  The precise choice of the coefficients is the subject of the next Section~\ref{sec:model-params}.

In the numerical experiments, we discretize equation~\eqref{eq:1} in space. We introduce a spatial grid for $\phi\in[0,\pi/2]$ with stepsize $\Delta\phi=\pi/(2I),I\in\mathbb{N}$. The gridpoints then are $\phi_{i}=i\Delta\phi,i=0,\ldots,I$. On these points, we denote the approximate  solution as $\bT_{i}(t)\approx T(\phi_{i},t)$. We then use standard finite differences, see Appendix~\ref{sec:discr}, to obtain a semi-discretization,
\begin{eqnarray}
\bT'(t)&=&f(t,\bT(t)),\quad t\ge0\label{eq:micro},
\end{eqnarray}
with $\bT(t),f(t,\bT(t))\in \mathbb{R}^{I-1}$ and some initial value $\bT(0)\in \mathbb{R}^{I-1}$. Equation~\eqref{eq:micro} is subsequently discretized in time, either by a forward Euler method with fixed time step or by a variable-step and variable-order linear multistep method, see Section~\ref{sec:numexp}.
 
\subsection{Spatial and temperature-dependent modeling\label{sec:model-params}}

 In the 1-D model, the parameters solar radiation $Q$ and albedo $\alpha$ can depend on colatitude $\phi$. The heat capacity $C$, albedo $\alpha$ and emissivity $\epsilon$ can depend on temperature, leading to nonlinearity in the models.  We first discuss the heat capacity, diffusion and emissivity (Section~\ref{sec:hde}), after which we go into more depth into the albedo (Section~\ref{sec:alb}) and the solar radiation (Section~\ref{sec:rad}).

\subsubsection{Heat capacity, diffusion and emissivity\label{sec:hde}} 

\begin{figure}
	\centering
	\includegraphics[scale=0.475]{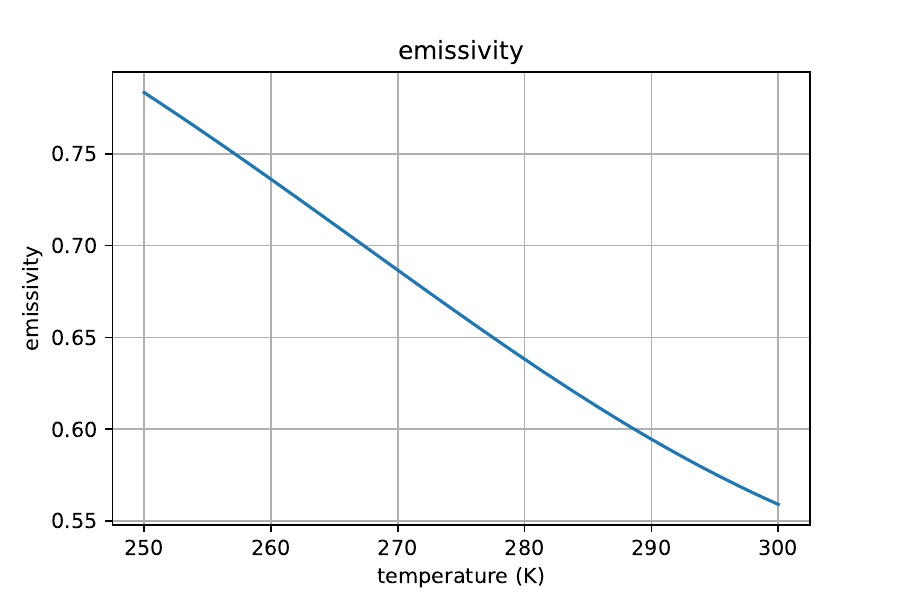}
	\includegraphics[scale=0.475]{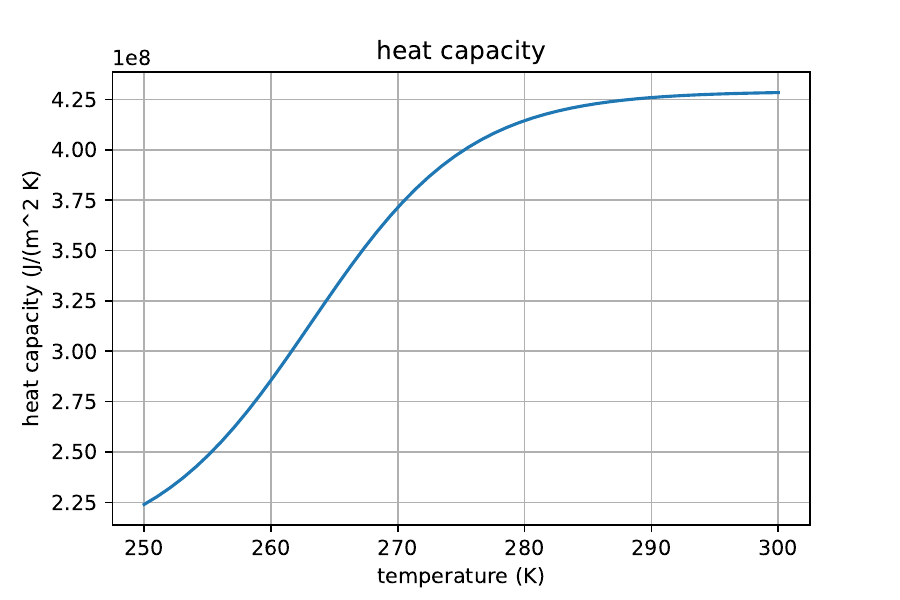}
	\caption{Left: Emissivity  $\epsilon=\epsilon(T)$. Right: Heat capacity $C=C(T)$. \label{fig:emmi}}
\end{figure}

\paragraph{Heat capacity.} For the heat capacity $C(T)$, we use  
the temperature-dependent function  from \cite[Section 4.5]{Bay91}:
\begin{eqnarray*}
C(T)&=&\left(C_{1}+C_{2}\tanh\left(C_{3}(T-T_{s})\right)\right)\times 10^{8}\,{\rm Jm}^{-2}{\rm K}^{-1}\end{eqnarray*}
with $C_{1}=3.14,C_{2}=1.15,C_{3}=0.08$, and $T_{s}=263.15$. 
{This nonlinear coefficient function is usually only considered in the 1-D model.}
It is shown in Figure \ref{fig:emmi}.

\paragraph{Diffusion coefficient.} For the diffusion coefficient $k(\phi,T)$, Ghil \cite[Eqn. (2e)]{Ghi75} used 
the following nonlinear function of colatitude and temperature 
\begin{eqnarray*}
k(\phi,T)&=&k_{1}(\phi)+k_{2}(\phi)g(T),\quad g(T)=\frac{c_{4}}{T^{2}}\exp\left(-\frac{c_{5}}{T}\right),
\end{eqnarray*}
with 
 given  data for the coefficient functions  $k_{1},k_{2}$ at the grid points (see \cite{Ghi75}) and  parameters
\begin{eqnarray*}
c_{4}
&=&6.105\times0.75\times\exp(19.6)\times10^{2}{\rm N  K  m}^{-2}
 \approx1.4891\times 10^{11} {\rm NKm}^{-2},\\
c_{5}&=&5.35\times10^{3}\,{\rm K}.
\end{eqnarray*}

Because this choice (especially the coefficients suggested in \cite{Ghi75}) does not guarantee positivity of $k$, we use a constant value, namely $k=0.591\,{\rm W m}^{-2}{\rm K}^{-1}$, taken from \cite[Section 4.5]{Bay91}. 
Note that, also in the literature cited above, mostly only one source of nonlinearity (either $C$ or $k$) has been chosen.

\paragraph{Emissivity.} For the emissivity $\epsilon$, we use
\begin{eqnarray*}
\epsilon(T)&=&1-\epsilon_{1}\tanh\left(\left(\frac{T}{T_{\epsilon}}\right)^{6}\right)
\end{eqnarray*}
with $\epsilon_{1}  = 0.5, T_{\epsilon}=284.15\,{\rm K}$, which  is  Ghil's   \cite{Ghi75} suggestion.

\subsubsection{Albedo\label{sec:alb}}
\begin{figure}
	\centering
	\includegraphics[scale=0.525]{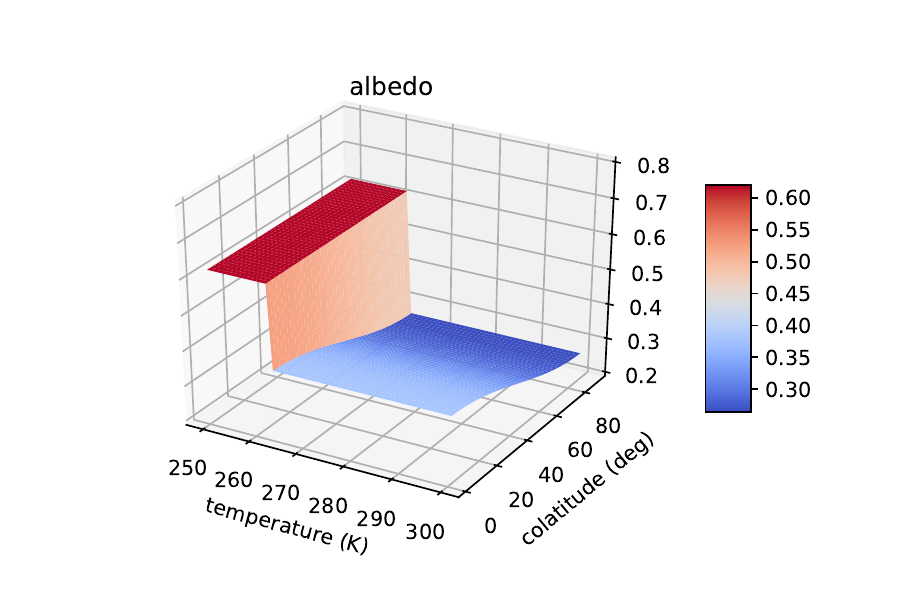}
	\includegraphics[scale=0.45]{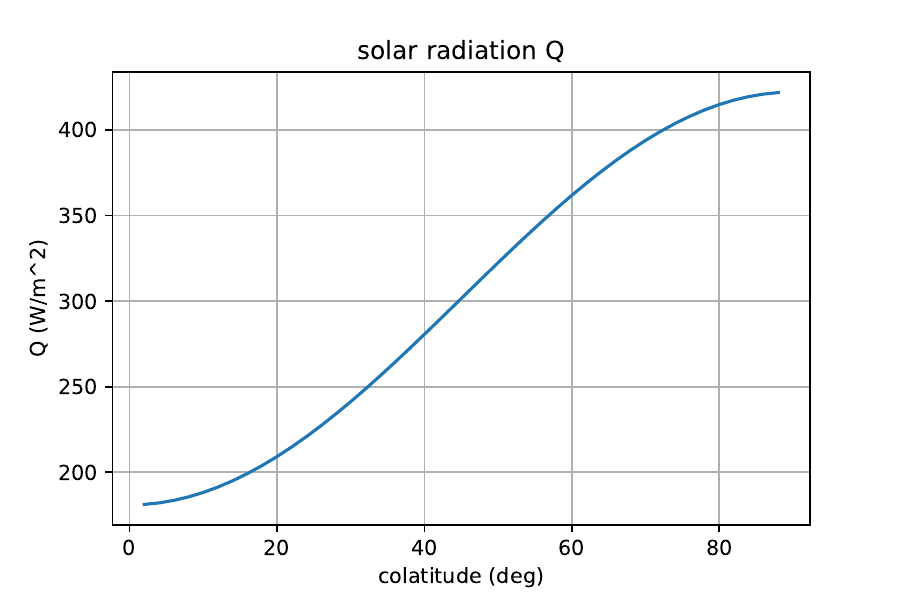}
	\caption{Left: Albedo  $\alpha=\alpha(\phi,T)$. Right: Solar radiation $Q=Q(\phi)$. \label{fig:albedo}}
\end{figure}

The albedo $\alpha(\phi,T)$ describes the fraction of solar radiation that is absorbed by the earth. Since this fraction depends on features such as clouds or snow and ice on the surface, the temperature-dependence is obvious. 
Piecewise constant and linear functions as well as smoothed counterparts are used in the literature. For the 1-D model,  spatial dependency reflects the bigger area covered by ice or snow (with higher albedo) near the poles. We use
 the model of North \cite{Nor75}, see also \cite[Section 4.5]{Bay91}:
\begin{eqnarray}
\label{eq:albedo}
\alpha(\phi,T)&=&
\left\{
\begin{array}{ll}
\alpha_{\rm max},& T\le T_{s},\\
\alpha_{1}+\alpha_{2}(\alpha_{3}\cos^{2}(\phi)-\alpha_{4}),&T> T_{s}, 
\end{array}
\right.
\end{eqnarray}
with $\alpha_{\rm max}=0.62,\alpha_{1}=0.303,\alpha_{2}=0.0779,\alpha_{3}=1.5,\alpha_{4}=0.5,$ and $T_{s}=263.15${K}.  {Here, $T_{s}$ is the temperature at the slowline. It can vary between approximately $-10$ and $0$ degree Celsius, see \cite[Section 3.2.2]{McgHen14}. The chosen value of $T_{s}$ corresponds to -10 degree Celsius.
The left-hand picture in} Figure \ref{fig:albedo} shows the resulting {discontinuous} function.

In the 0-D version, we restrict to a simple step function
\begin{eqnarray}
\label{eq:albedo-step}
\alpha(T)&=&\left\{
\begin{array}{ll}
\alpha_{\rm max},&T\le T_{s} ,\\[0.2em]
\alpha_{\rm min},\quad&T> T_{s},
\end{array}\right.
\end{eqnarray}
which is also mentioned in  \cite[Section 10.2]{GhiChi87} 
and \cite[(3.10)]{McgHen14}. 
{This step function is a simplification, assuming an instantaneous change in albedo  due to ice and snow melting at temperature $T_{s}$.} The parameter $\alpha_{\rm min}$ in the 0-D model is not present in the 1-D model. We chose $\alpha_{\rm min}=0.275$, which gives a good approximation of the long-time steady state obtained by the micro model.

\subsubsection{Solar radiation\label{sec:rad}}

\paragraph{Space dependence.} Clearly, the incoming solar radiation $Q(\phi)$ depends on the colatitude of the considered position on the Earth's surface. 
We use a second order polynomial in the variable $x=\sin\phi$, see \cite[Section 4.3]{Sto14} or \cite[Section 4.2]{Bay91}),
 \begin{eqnarray}
 \label{eq:Q}
Q(\phi) &=& \frac{S}{4}\left(Q_{1} + Q_{2} \sin^2\phi\right)
\end{eqnarray}
with $Q_{1}=0.5294$ and $Q_{2}=0.706$, see Figure \ref{fig:albedo}.
Ghil \cite{Ghi75} took data of which the spatial distribution is similar to the one of the polynomial above. 

 For the 0-D model, we use the average of~\eqref{eq:Q} over the half-sphere.
The spatial mean of a quantity depending on colatitude $\phi$ is obtained by integrating it over $\phi\in[0,\pi/2]$ and dividing by the measure of the area. In spherical coordinates the arc length of the circle at constant colatitude $\phi$ is given by $\pi r\sin\phi$ ($r$ being again the Earth's radius). This gives
\begin{eqnarray}\label{eq:spatial_mean}
\overline{Q}&=&\left(\int_{0}^{\frac\pi2}Q(\phi)\sin\phi \,d\phi\right)\left(\int_{0}^{\frac\pi2}\sin\phi \,d\phi\right)^{-1}
=\int_{0}^{\frac\pi2}Q(\phi)\sin\phi \,d\phi,
\end{eqnarray}
since the integral in the denominator equals 1. The integral for $Q$ given by the formula \eqref{eq:Q} can be evaluated exactly.  
Since we also made experiments with the date provided by Ghil, we use a numerical approximation by the trapezoidal quadrature rule on the numerical grid on which equation~\eqref{eq:1} is solved.

\begin{figure}
	\centering
	\includegraphics[scale=0.5]{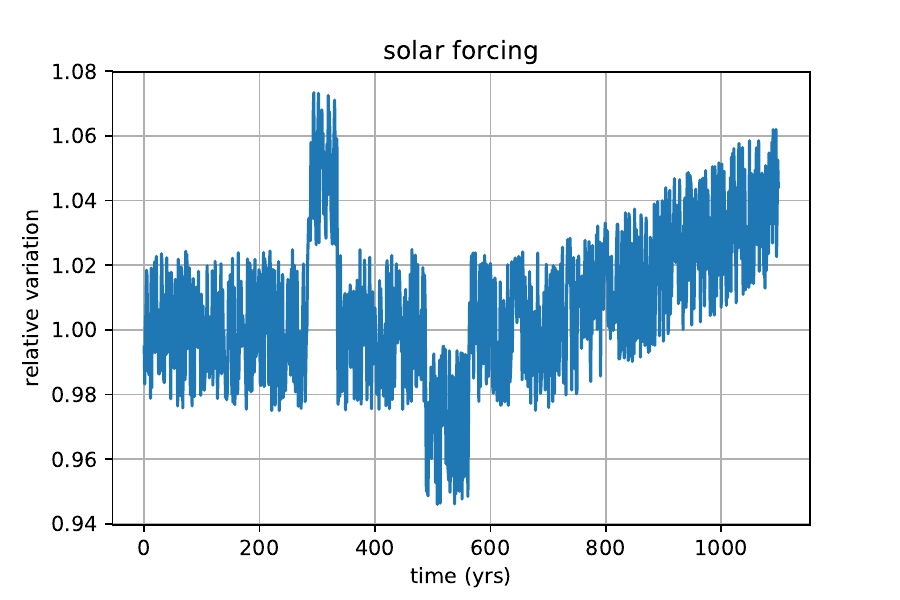}
	\caption{Introduced variation of solar forcing with four jumps, a linear trend after year $700$, and a random variation in the whole time interval $[0,1000]$, compare  \eqref{Qdist0}, \eqref{Qdist1} and \eqref{Qdist2}. \label{fig:fluct-1}}
\end{figure}

\paragraph{Temporal variation.} 
Temperature data reconstructions over one or more glacial cycles typically show relatively sharp gradients which are superimposed by small-scale fluctuations \cite{Brook2008}. Note that, in this context,  ``sharp'' has to be interpreted in relation to  the considered time ranges of  several hundred thousands of years. Furthermore,  a global warming trend 
can be seen in temperature observations in the last  decades. 
Our motivation for the use of the parareal method is to capture such multi-scale behavior in climate  models. Thus we 
additionally varied $Q$ in the 1-D model by adding some jumps, a linear trend and small-scale random fluctuations, described by a time-dependent function $\Delta Q$. We write
\begin{align}
\label{Qdist0}
Q_{\Delta}(\phi,t)&=Q(\phi)(1 + \Delta Q(t)),
\end{align}
in which the perturbation $\Delta Q(t)$ has the following form:
\begin{align}
\label{Qdist1}
\Delta Q(t)&=\sum_{i=1}^{2}q_{i}\chi_{[a_{i},b_{i}]}(t)+q_{3}\chi_{[a_{3},b_{3}]}(t)(t-a_{3})+q_{4}R(t).
\end{align}
Here $\chi_{[a,b]}$ is the characteristic function on the interval $[a,b]$ and $\{R(t):t\in\{1,\ldots,t_{\rm end}\}\}$ is a set of uniformly distributed random numbers in $[-1,1]$. These were computed once and  then fixed for all experiments. For our tests, we {used}  the values
\begin{align}
\label{Qdist2}
\begin{array}{llllll}
q_{1}=0.05,&a_{1}=283,&b_{1}=335,&
q_{2}=- 0.03,&a_{2}=487,&b_{2}=564,\\
q_{3}=0.0001,&a_{3}=700,&b_{3}=1000,&{q_{4}=0.025}.
\end{array}
\end{align}
In the time integration schemes, the values of $R(t)$ were interpolated linearly  whenever needed at non-integer time instants. Figure \ref{fig:fluct-1} shows the relative deviation from the constant incoming solar radiation $Q$ as function of time. {The random fluctations can be seen as effect of sunspot activity, with amplified magnitude since the relative variation due to the latter is in the magnitude of $\approx 10^{-3}$.}

For the 0-D model, we tested two choices to handle the time-dependence of $Q_{\Delta}(t)$: 
\begin{enumerate}[(1)]
\item a spatially averaged value that takes into account the time-dependence of $Q_{\Delta}(\phi,t)$ (ignoring the random fluctuations); 
\[
\overline{Q}_{\Delta}(t)= \int_{0}^{\frac\pi2}Q_{\Delta}(\phi,t)\sin\phi \,d\phi
\]
\item a constant value $\overline{Q}_{\Delta}(t)=\overline{Q}$ for all $t$. 
\end{enumerate}

\subsection{Existence and uniqueness of solutions\label{sec:model-exist}}

\begin{figure}
	\centering
	\includegraphics[scale=0.475]{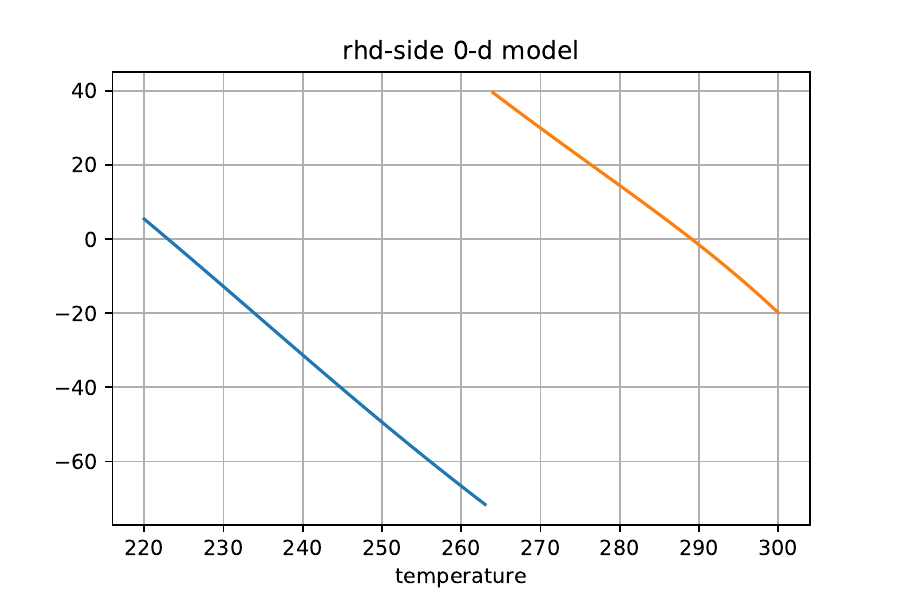}
	\includegraphics[scale=0.475]{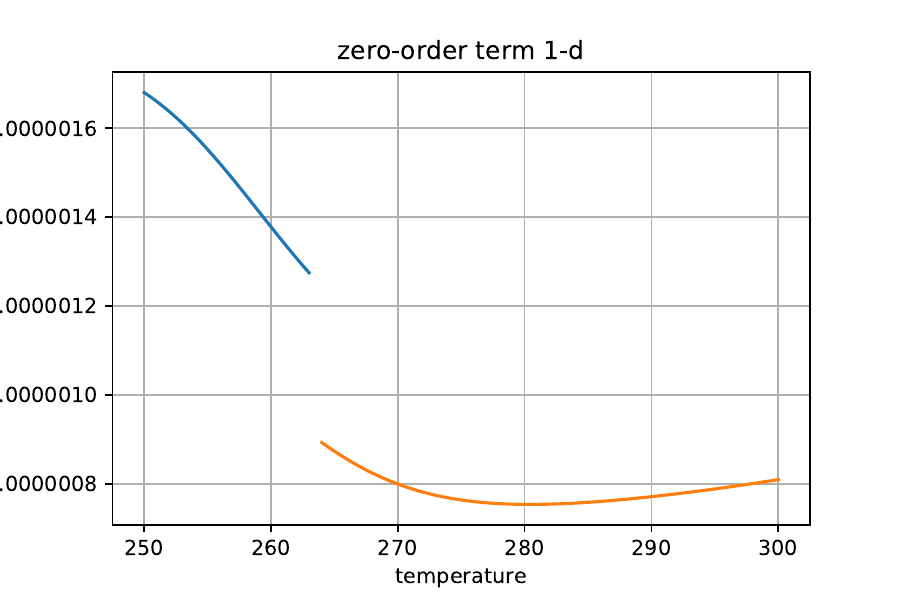}
	\caption{Left: Numerator of the right-hand side of the 0-D model. Right: Zero-order term of the 1-D model, here for spatially averaged values of $\alpha,Q$. \label{fig:ana}}
\end{figure}

In this section, we briefly summarize existence and uniqueness results of the two models. 
With the temperature-dependent modeling introduced above, the 0-D model takes the form
\begin{eqnarray}
\label{eq:0-D2}
T'=
\frac{(1-\alpha(T))\overline{Q}
 - \epsilon(T)\sigma T^4}{C(T)}.
\end{eqnarray}
With the albedo being the step-function \eqref{eq:albedo-step}, the right-hand side does not fulfill the classical assumptions of the theorems of Peano or Picard-Lindel\"of for existence and uniqueness. On the other hand, we have $\tanh(x)\in[-1,1], x\in\mathbb{R}$, and thus 
$C(T)>0$ for $ T\in\mathbb{R}_{\ge0}.$ Hence, the sign of $T'$ can be deduced from the  numerator of the right-hand side of \eqref{eq:0-D2}, compare Figure \ref{fig:ana}: For $T<T_{s}$ we have $T'<0$ and for $T>T_{s}$ we have $T'>0$. Thus any solution of the 0-D initial value problem for \eqref{eq:0-D2} with $T(0)=T_{0}$ will always remain in either $(0,T_{s})$ or $(T_{s},\infty)$, depending in which interval $T_{0}$ lies. In each of these intervals,  the right-hand side of the equation is continuously differentiable  w.r.t. $T$ and thus locally Lipschitz continuous. Thus the Picard-Lindel\"of  Theorem is applicable and gives existence and uniqueness of the solution. 

The 1-D model is a semilinear parabolic PDE of the general form
$$
C(T)\,T_{t}+\left(kT_{\phi}\right)_{\phi}+g(\phi,T)=f,
$$
where subscripts denote partial derivatives.
Since $C$ is positive it can be re-written as
$$
T_{t}+\frac{1}{C(T)}\left(kT_{\phi}\right)_{\phi}+\tilde g(\phi,T)=0.$$
 It is not the aim of this paper to perform an analysis of the equation. We just  note that the zero-order term
$$
\tilde g(\phi,T)=\frac{(\alpha(\phi,T)-1)Q(\phi)
 +\epsilon(T)\sigma T^4}{C(T)}
 $$
 shows -- besides the discontinuity at $T=T_{s}$ coming from the albedo function -- a non-monotonicity w.r.t. $T$ for $T\ge T_{s}$, compare Figure \ref{fig:ana}. Classical existence theorems (see e.g. \cite[Section 9.2, Theorem 2]{Eva98}, \cite[Section 7.3, Lemma 5.3]{Tro10}) require either Lipschitz continuity or monotonicity of $\tilde g$. This lack of theoretical results shows the typical analytical properties of nonlinear climate models, even if they are low-dimensional as in this case.

\section{Micro/macro parareal algorithm\label{sec:parareal}}

In this section, we propose the micro/macro parareal algorithm that will be studied further on.  It is a straightforward generalization of the parareal algorithm proposed in \cite{LiMaGa01}; see also \cite{LeLeSa13}, where a similar algorithm was proposed in the context of singularly perturbed ordinary differential equations.  We first introduce the necessary notation (Section~\ref{sec:parareal-notation}), after which we outline the algorithm (Section~\ref{sec:parareal-alg}).

\subsection{Notations\label{sec:parareal-notation}}

We introduce a time discretization $(t_n)_{n=0}^{N}$, with $t_n=n\Delta t$, as well as the numerical approximations $\bT^n\approx \bT(\phi,t_n)$ of the 1-D model \eqref{eq:1} (or, more precisely, its spatial discretization~\eqref{eq:micro}), and $T^n\approx T(t_n)$ of the macroscopic model \eqref{eq:macro}, respectively. 

\paragraph{Fine-scale and macro propagators.}
The micro/macro parareal algorithm makes use of a \emph{micro-scale propagator}, that advances the microscopic model \eqref{eq:micro} over a time-step $\Delta t$, 
\begin{equation}\label{eq:micro-propagator}
\bT^{n+1}=F_{\Delta t}(\bT^n),
\end{equation}
and, similarly, a \emph{macro propagator} for the macroscopic model
\eqref{eq:macro},
\begin{equation}\label{eq:macro-propagator}
T^{n+1}=G_{\Delta t}(T^n).
\end{equation}
For now, we consider a forward Euler time discretization with time step $\Delta t$ for the macro propagator and with time step $\delta t \ll \Delta t$ for the micro-scale propagator.

\paragraph{Lifting, restriction and projection.} Furthermore, we introduce operators that connect the microscopic and macroscopic levels of description. The \emph{restriction} operator
\begin{equation}
\R:  T \mapsto T = \R(\bT),
\end{equation} 
maps a microscopic state to the corresponding macroscopic state.
For the model problem \eqref{eq:macro}-\eqref{eq:micro}, this restriction operator simply becomes the averaging formula~\eqref{eq:spatial_mean}.
For notational convenience, we also introduce the complement of the restriction operator, 
 \[
 \R^\perp(\bT):= \dfrac{\bT}{\R(\bT)},
 \]
such that we can write $T(\phi,t)=\R(T(\phi,t))\cdot\R^\perp(T(\phi,t))$. 

Conversely, to reconstruct a microscopic state from a given macroscopic state, we distinguish between a \emph{lifting} operator, and a \emph{projection} operator.  A lifting operator $\L$ initializes a microscopic temperature profile $T(\phi,t)$ that is uniquely determined by the given spatially averaged temperature $T(t)$,
\begin{equation}
\L: T \mapsto \bT = \L(T).
\end{equation}
For the model problem \eqref{eq:micro}-\eqref{eq:macro}, for instance, we could choose
\begin{equation}\label{eq:lift}
	T(\phi,t)=\L(T(t)):=T(t)\Psi(\phi),
\end{equation}
with $\Psi(\phi)$ an arbitrary function such that $\R(\Psi(\phi))=1$. 
Clearly, one requires $\R \circ \L=\Id$.  

In contrast, one may also \emph{match} a ``nearby'' (\emph{prior}) temperature profile $\bT^*(\phi)$ with a desired macroscopic spatially averaged temperature $T$. Then, the result is not uniquely determined by the macroscopic state $T$, but depends also on the prior $\bT^*(\phi)$. We call the resulting reconstruction operator a matching operator,
\begin{equation}
\P: T,\bT^*(\phi) \mapsto \bT = \mathcal{P}(T,\bT^*(\phi)).
\end{equation}
This operator projects a microscopic temperature profile $\bT^*(\phi)$ onto the manifold of microscopic temperature profiles consistent with the macroscopic spatially averaged temperature $T$. 
Here, we require $ T = \left(\R \circ \P\right)(T,\bT^*(\phi))$, for any $\bT^*(\phi)$.  
Additionally, a so-called \emph{self-consistency} property is of particular importance.
\begin{defn}[Self-consistency] \label{def:self-consistent} A projection operator $\P : T,\bT^*(\phi) \mapsto \bT(\phi) = \mathcal{P}(T,\bT^*(\phi))$, is called self-consistent if, and only if,
\begin{equation}
\forall \bT(\phi): \R(\bT(\phi))=T \;\; \Rightarrow \;\; \bT(\phi)=\P(T,\bT(\phi)).
\end{equation}
\end{defn}
 When this property holds, a microscopic temperature profile is not altered if it is projected onto a macroscopic spatially averaged temperature with which it is already consistent.
As a guideline, the matching should be such that $\mathcal{R}(\bT(\phi))=T$, while requiring $\bT(\phi)$ to be as close to $\bT^*(\phi)$ as possible, in a sense to be made precise for the problem at hand.   
For the model problem \eqref{eq:micro}-\eqref{eq:macro}, we choose the matching operator as
\begin{equation}\label{eq:proj}
	\bT(\phi)=\P(T,\bT^*(\phi)):=T\dfrac{\bT^*(\phi)}{\R(\bT^*(\phi))}=T\cdot \R^\perp(\bT^*(\phi)).
\end{equation}

\subsection{Algorithm\label{sec:parareal-alg}}

The parareal algorithm iteratively constructs approximations on the whole time domain.  We denote by $T_k^n$ the approximate solution at time $t_n$, obtained during the $k$-th parareal iteration. 
We start from an initial condition $\bT(\phi,0)=T_0$, and create an initial approximation on the whole time interval by using the macro propagator, i.e. we restrict $T_0 = \R(\bT_0)$, and compute
\begin{equation}
T_0^{n+1}=G_{\Delta t}(T_0^n), \qquad T^0_0 = T_0.
\end{equation}
We then lift this macro approximation to the micro scale,
\begin{equation}
\label{eq:lift1}
\bT_0^n = \L(T_0^n).
\end{equation} 
We now have an initial approximation of the microscopic solution at each of the time instants $t_n$, $1\le n \le N$.
The parareal iterations then proceed as follows.
\begin{enumerate}[a)]
\item Compute (in parallel) the time propagation at each time instance, using the macro and micro-scale propagators,
\begin{align}
\bar{T}_k^{n+1}&=G_{\Delta t}(T_k^n), \\
\bar{\bT}_k^{n+1}&=F_{\Delta t}(\bT_k^n).\label{eq:micro_prop}
\end{align}
\item Compute the jumps (the difference between the two propagated values) \emph{at the macro level},
\begin{align}\label{eq:parareal-jumps}
J_k^{n+1}= \R(\bar{\bT}^{n+1}_k)-\bar{T}^{n+1}_k.
\end{align}
\item Propagate (serially) the macro jumps over the whole time domain using the macro propagator,
\begin{align}\label{eq:parareal-propagate}
T_{k+1}^{n+1} &= G_{\Delta t}(T^n_{k+1})+J^{n+1}_k,
\end{align}
and reconstruct the micro-scale state from the corrected macro state,
\begin{align}\label{eq:parareal-projection}
\bT_{k+1}^{n+1}=\P(T_{k+1}^{n+1},\bar{\bT}_k^{n+1}).
\end{align}
\end{enumerate}
Remark that the micro-scale state is reconstructed by projecting the intermediate value $\bar{\bT}_k^{n+1}$ onto the corrected macro value $T_{k+1}^{n+1}$, using the matching operator~\eqref{eq:proj}.

\section{Setup of numerical experiments}
\label{sec:numexp}
In this section, we describe the numerical experiments that we performed.  In Section~\ref{sec:time-choices}, we discuss the two options for the time integration of the micro and macro models.  In Section~\ref{sec:macro-choices}, we discuss the two versions of the macro model that will be considered. 
The source code is available at \url{https://doi.org/10.5281/zenodo.1287561}.

\subsection{Time integration\label{sec:time-choices}}
\label{sec:time-integration}
\paragraph{Two time discretization methods.} We simulate up to time $t_{\rm end}=1000$, using an initial value of $285\,{\rm K}$ for all computations, in 0-D as well as uniform in space for the 1-D model.
For the micro-scale propagator $F_{\Delta t}$, we use the spatial discretization~\eqref{eq:micro} for the 1-D model, using $I=45$. For the time discretization in the micro 1-D model, we use the built-in \verb#odeint# routine in Python's \verb$scipy$ library, which is the routine \verb$lsoda$ \cite{Pet83} from the \verb$odepack$ library \cite{Hin83}. This method automatically switches between an Adams-Bashforth method (for non-stiff problems) or a backward Differentiation Formula (BDF, for stiff ones). The 1-D model is stiff, thus nearly in all cases and at all time instants the implicit version was chosen by the algorithm. 
The \verb#odeint# routine also uses a variable time step, based on a user-prescribed error tolerance. We use $ \varepsilon_{\rm rel}= \varepsilon_{\rm abs}=10^{-6}$  for both relative and absolute tolerances. This means that a  time-step is accepted if the predicted error $e_{n}$ in time instant $t_{n}$ satisfies $$\|e_{n}\|_{\infty}\le \varepsilon_{\rm rel}\|\bT_{n}\|_{\infty}+\varepsilon_{\rm abs}.$$  

Note that, because of the time-step adaptation strategy, we have no control over the time steps taken during time integration. In particular, it is very unlikely that the same time step will be taken in subsequent parareal iterations (or in a serial micro simulation over the full time interval). As a consequence, the exactness property of the parareal method becomes very hard to check numerically.  We therefore also perform numerical experiments in which the micro model is integrated using the explicit Euler method with fixed stepsize. {The stability  condition for explicit time-stepping of a linear diffusion equation leads to an upper bound of $\frac{h^{2}}{2k}$ for the step-size in time, here with  $h=\frac{\pi}{2I}=\frac{\pi}{90}$ and $k$ being the diffusivity. In our setting, this results in a value of $\approx 10^{-3}$. This rough calculation neglects all nonlinearities. It turned out that a step-size of $5\times10^{-3}$ gave reasonable results for all experiments.}

\begin{figure}
	\centering
	\includegraphics[scale=0.475]{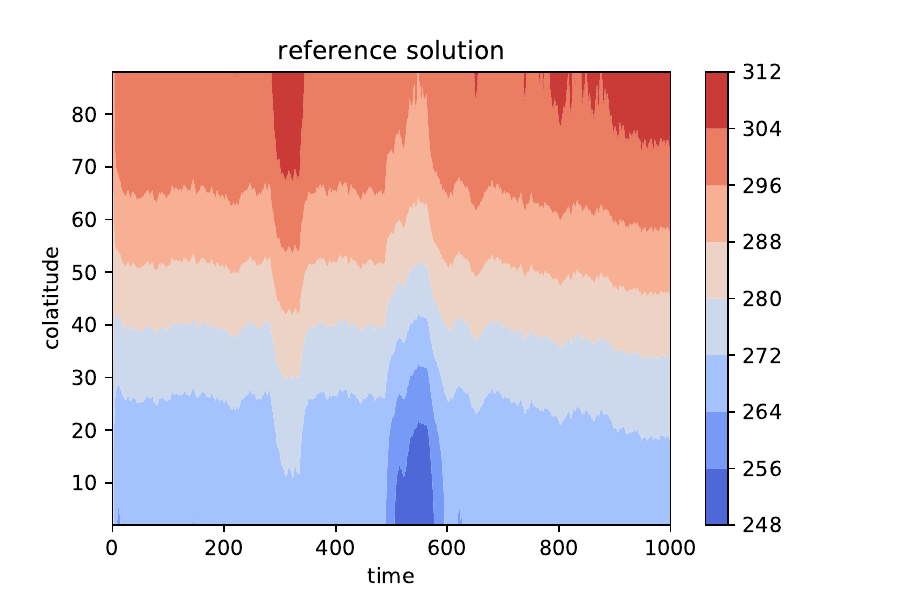}
        \includegraphics[scale=0.475]{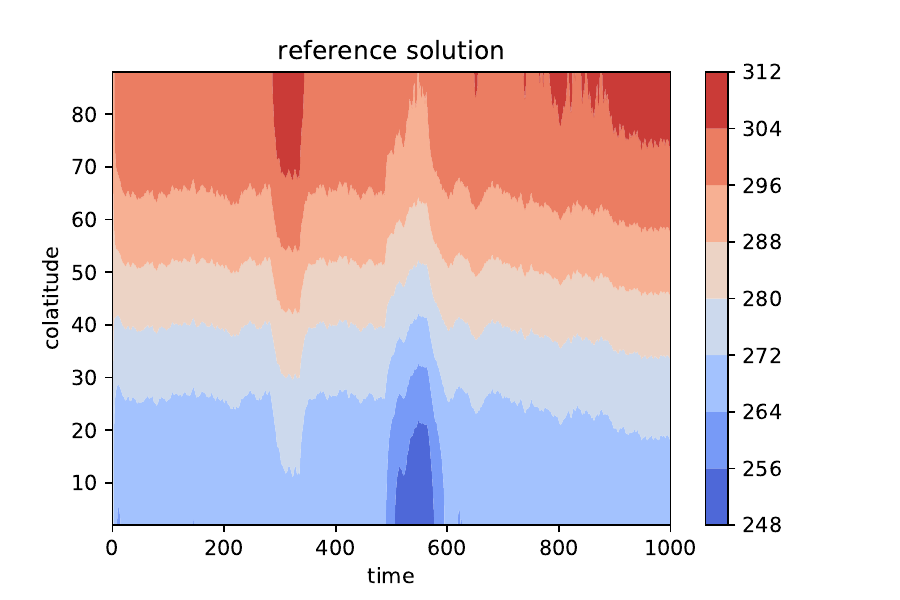}
	\caption{1-D macro reference solutions obtained by the Euler method with constant stepsize 0.005 (left) and the adaptive  {\tt lsoda} library routine. \label{fig:refsol}}
\end{figure}

\paragraph{Reference solutions.} We compute wo 1-D micro solutions sequentially using each of the above-described time discretization methods, and use them as reference solutions to compute the differences to the respective parareal solutions. These reference solutions are are denoted by ${\bf T}^{*}(t)$. and shown in figure \ref{fig:refsol}. 

\paragraph{Macro propagator.} As the macro-scale propagator $G_{\Delta t}$ for the 0-D model, we again used the explicit Euler method, but now with the constant stepsize  $10.0$.  

 \subsection{Two versions of macro model\label{sec:macro-choices}}
\begin{figure}
	\centering
	\includegraphics[scale=0.465]{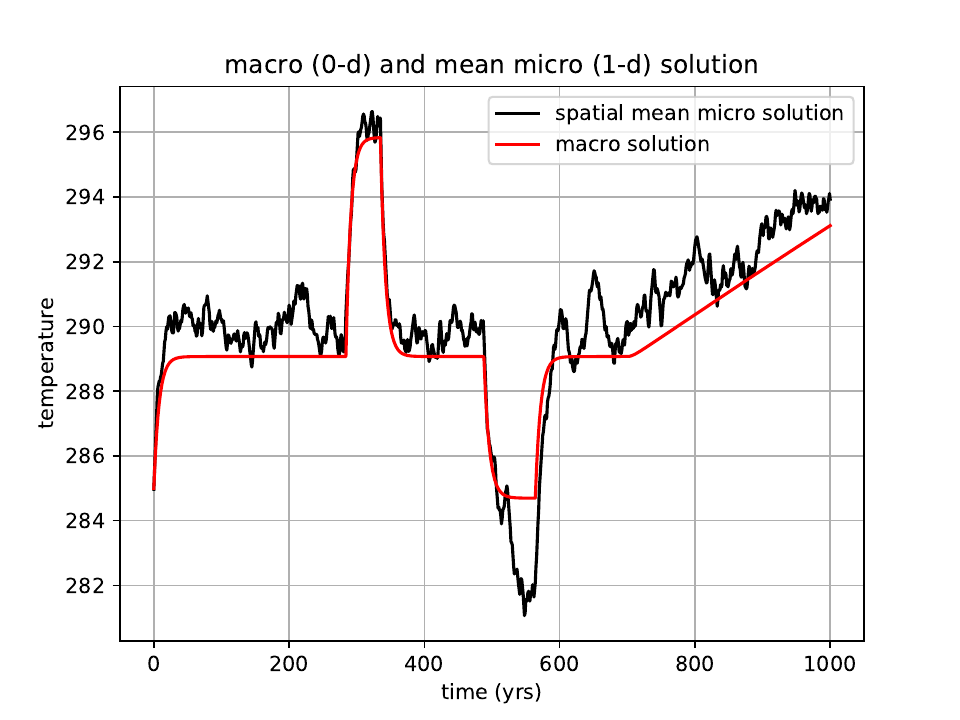}
		\includegraphics[scale=0.465]{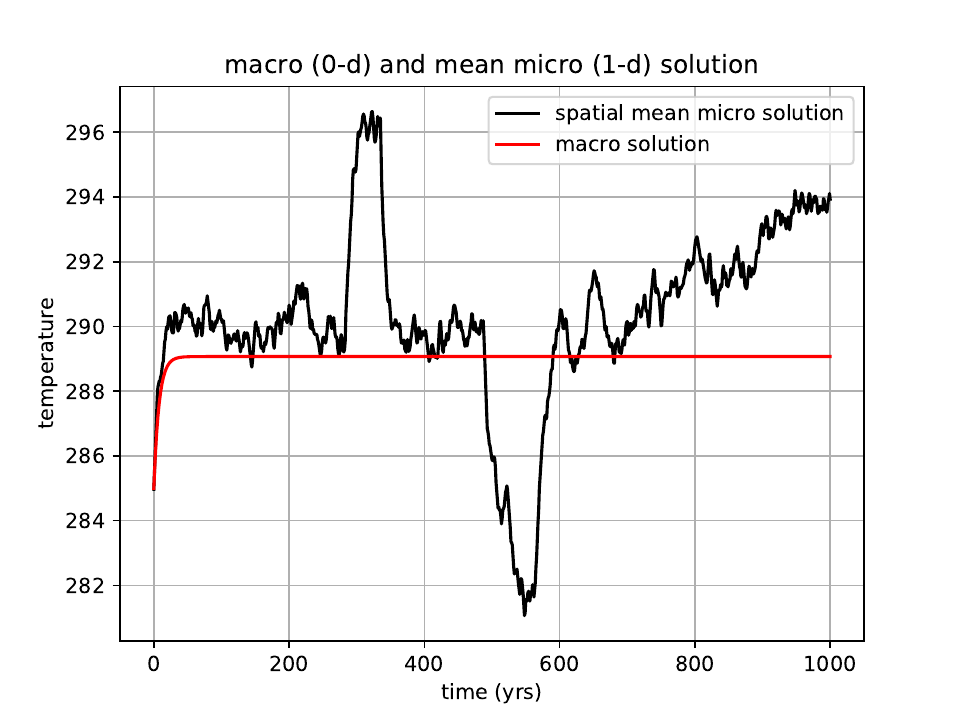}
	\caption{Spatially averaged 1-D micro solution $\overline{\bT}$ with introduced perturbation compared to macro solution $T$ with variable (left) and constant forcing (right). \label{fig:fluct-2}}
\end{figure}

Besides varying the microscopic time discretization, we also use two versions of the  0-D macro model. One has a similar temporally varying forcing as the micro model, but without the small-scale random fluctuations. This means we use  $\Delta Q$ as in \eqref{Qdist1} and \eqref{Qdist2}, but with $q_{4}=0$. The other version of the macro model has a constant solar forcing, i.e., $\Delta Q=0$. The  1-D micro reference solution (with varying forcing) obtained with the Euler method (the same as in Figure \ref{fig:refsol} on the left, now averaged in space) is compared in Figure \ref{fig:fluct-2} with these two macro solutions.  As can be seen and is obvious by construction, the macro model with constant solar forcing is an even more coarse and inaccurate approximation of the 1-D micro model. It does not follow neither the jumps nor the linear trend of the forcing (which is {caught} by the micro model), whereas the macro model with variable forcing does, at least to some extend. Our aim was to see how much this fact influences  the parareal convergence.

\subsection{Summary of numerical setup}

We thus end up with four configurations, namely using either Euler or an adaptive library routine as micro solver and using either the macro model with or without added perturbation $Q_{\Delta}$ in solar forcing.
For all these four settings, we varied the parareal time-step $\Delta t\in\{10,20,25,40,50,100\}$ (corresponding to $N\in\{100,50,40,25,20,10\}$) and compared the results w.r.t. difference of the parareal solution $\bT_{k}$ in iteration $k$ and the serially computed micro reference solution  $\bT^{*}$. The differences were always evaluated  at all integer time instants $t\in \{1,\ldots,t_{\rm end}\}$.    

With the tolerances set as in Section~\ref{sec:time-choices}, the library routine \verb#odeint# and the forward Euler method require approximately the same time for a computation. In this setting, the  0-D macro model with constant forcing requires (approximately and on average (over all performed runs) only a fraction $1/700$ of the micro computation time. The  macro model with temporally varying forcing takes about 2-3 times as much time as macro model with constant forcing. Recall that the 1-D micro model always uses the varying forcing.


\section{Numerical results\label{sec:results}}

\subsection{Results using Euler method as micro propagator}
\begin{figure}
	\centering
	\includegraphics[scale=0.3]{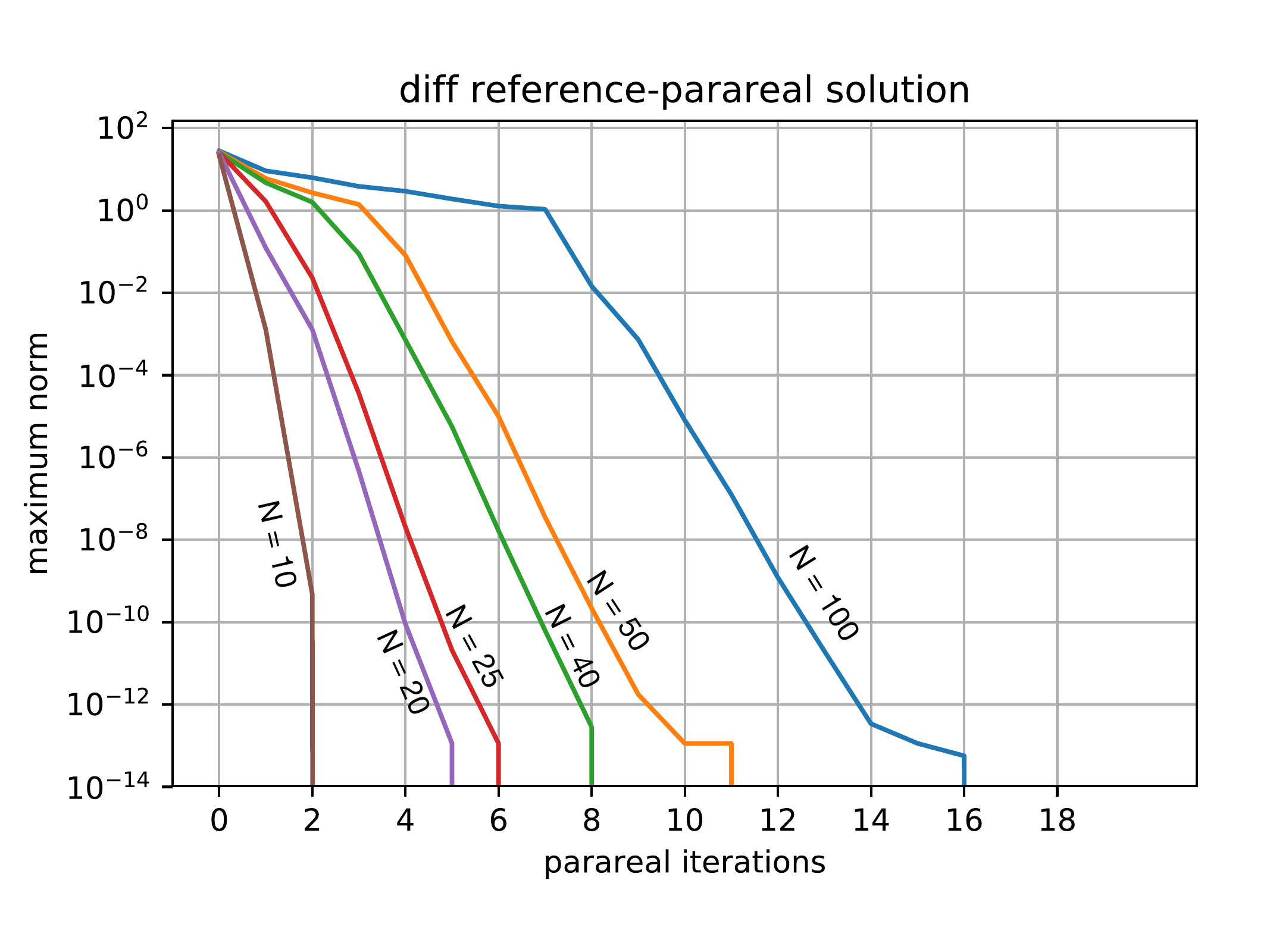}
	\caption{Convergence of parareal method using Euler method for the micro model and  the macro model with  temporally varying solar forcing (depicted in left plot of Figure \ref{fig:fluct-2}) for different numbers $N$ of parareal subintervals. Plotted are maximum norms over space and time.  \label{fig:conv-1}}
\end{figure}

\begin{figure}
	\centering
	\includegraphics[scale=0.3]{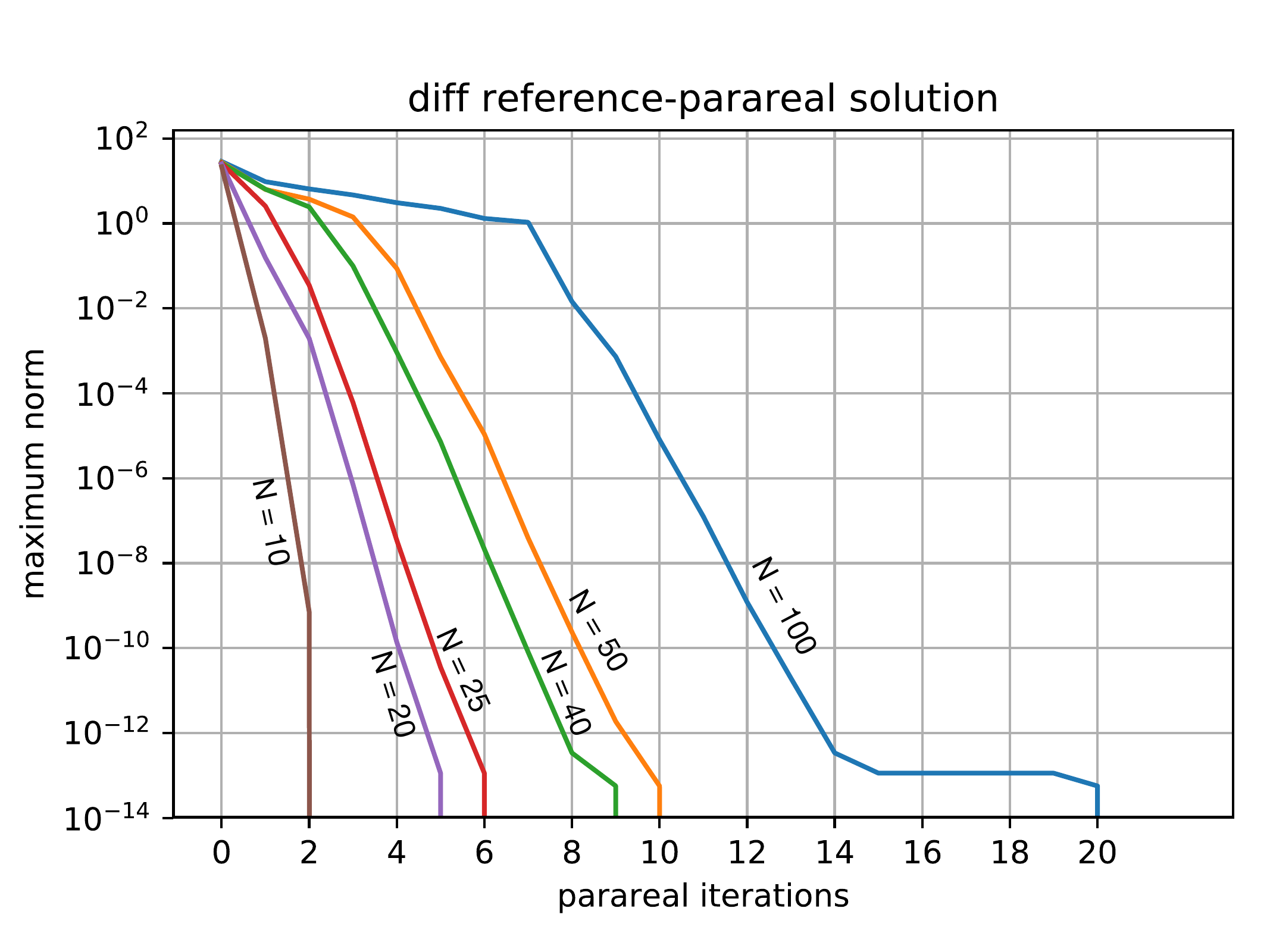}
	\caption{Same as Figure \ref{fig:conv-1}, but   using  the macro model with constant solar forcing (depicted in right plot of Figure \ref{fig:fluct-2}).\label{fig:conv-2}}
\end{figure}

Figures \ref{fig:conv-1} and \ref{fig:conv-2} show the convergence of the parareal method to the micro reference solution using the Euler method  for the micro model. 
It can be seen that the differences reach machine precision, see  also Table \ref{table:1}. A reasonable difference level of $10^{-2}$ (talking about temperature which is the range of 200-300 K)  is reached even faster.

To investigate the theoretical gain in computational effort, we assume that  $N$ processors are used and discard the effort for communication and macro propagator. It can be seen from the values in the table that  a  reduction by  a factor of approximately  10    is reached. This refers to a computation up to a difference of the given reasonable tolerance of $10^{-2}$ {that was achieved by the settings of the tolerances and step-sizes mentioned in Section \ref{sec:time-integration}.} For computations up to machine precision the factor is approximately 3 to 5. 

\begin{table}
\centering
\begin{tabular}{|c|c|r|r|r|r|r|r|}
\hline
Micro&$N$ (number of parareal subintervals)&100&50&40&25&20&10\\
\cline{2-8}
propagator&$\Delta t$ (length of parareal subintervals)&10&20&25&40&50&100\\
\hline
explicit&\multicolumn{7}{|c|}{macro model: temporally varying forcing}\\
\cline{2-8}
Euler&$\min \{k: \|\bT_{k}-\bT^*\|_\infty<eps\}$&17&12&9&7&6&3\\\cline{2-8}
&$\min \{k: \|\bT_{k}-\bT^*\|_\infty<10^{-2}\}$&9&5&4&3&2&1\\\cline{2-8}
&\multicolumn{7}{|c|}{macro model: constant forcing}\\\cline{2-8}
&$\min \{k: \|\bT_{k}-\bT^*\|_\infty<eps\}$&21&11&10&7&6&3\\\cline{2-8}
&$\min \{k: \|\bT_{k}-\bT^*\|_\infty<10^{-2}\}$&9&5&4&3&2&1\\\hline
{\tt lsoda}&\multicolumn{7}{|c|}{macro model: temporally varying forcing}\\\cline{2-8}
library&$\min \{k: \|\bT_{k}-\bT^*\|_\infty<10^{-2}\}$&12&5&4&3&3&1\\\cline{2-8}
routine&\multicolumn{7}{|c|}{macro model: constant forcing}\\\cline{2-8}
&$\min \{k: \|\bT_{k}-\bT^*\|_\infty<10^{-2}\}$&10&5&4&3&3&1\\\hline
\end{tabular}
\caption{Convergence of the micro/macro parareal method to a reasonable tolerance of $10^{-2}$, using the Euler method and the adaptive library method in the micro model. For the Euler method also the convergence  to machine precision (in double precision IEEE arithmetic, $eps\approx 2.2\times 10^{-16}$) is shown.  $\bT_{k}$ denotes the parareal solution in the $k$-th iteration, $\bT^{*}$ the reference solution.
\label{table:1}}
\end{table}

\begin{figure}
	\centering
	\includegraphics[scale=0.465]{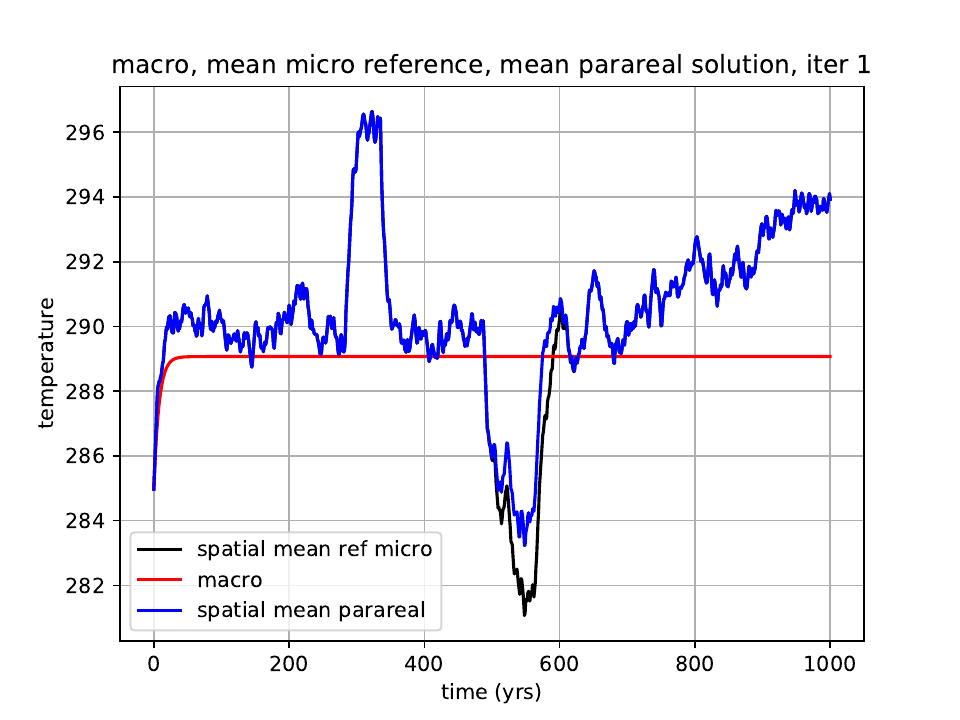}
	\includegraphics[scale=0.465]{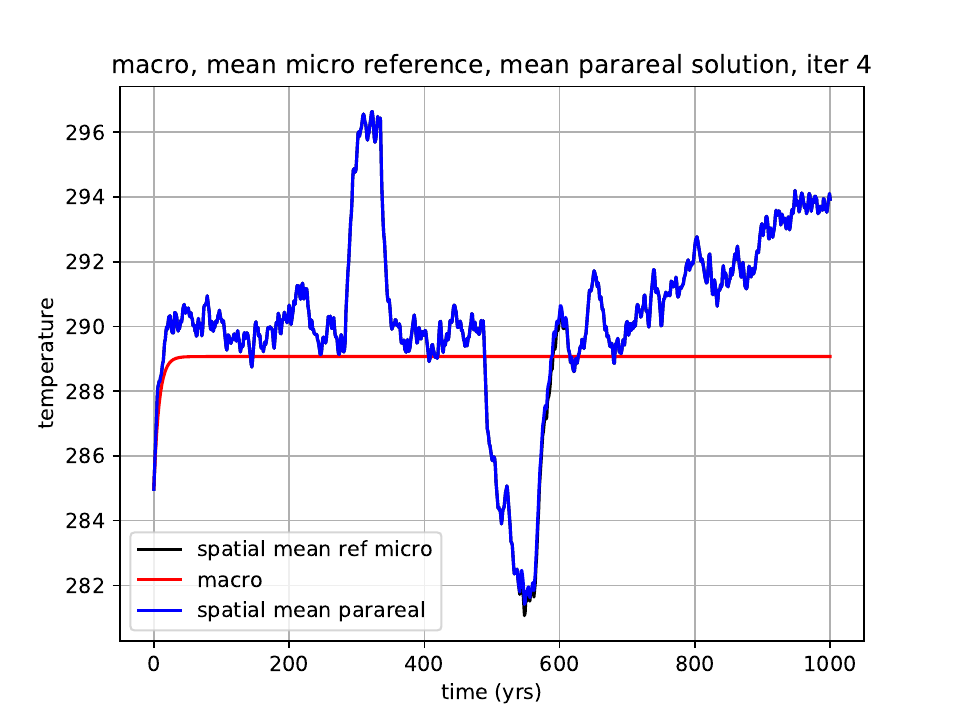}
	\includegraphics[scale=0.475]{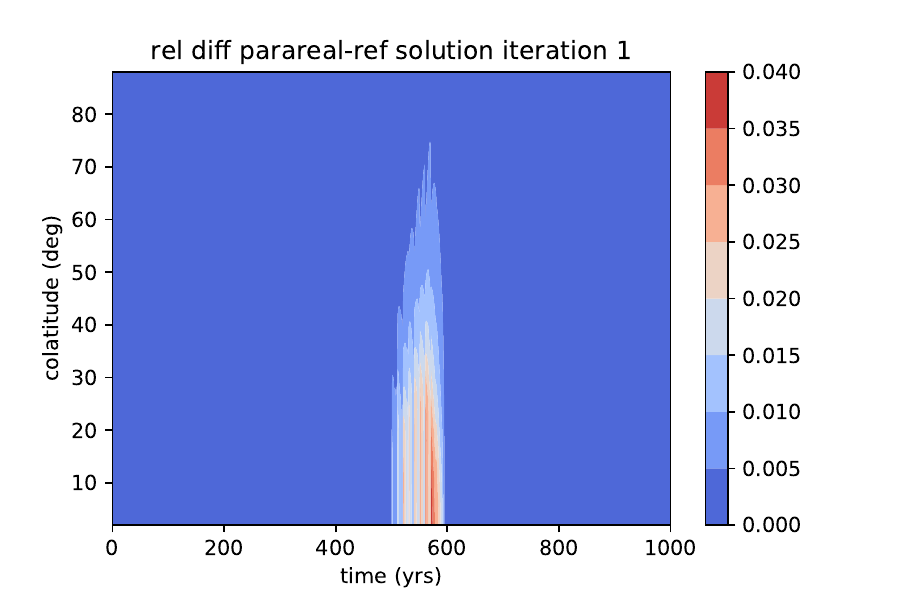}
	\includegraphics[scale=0.475]{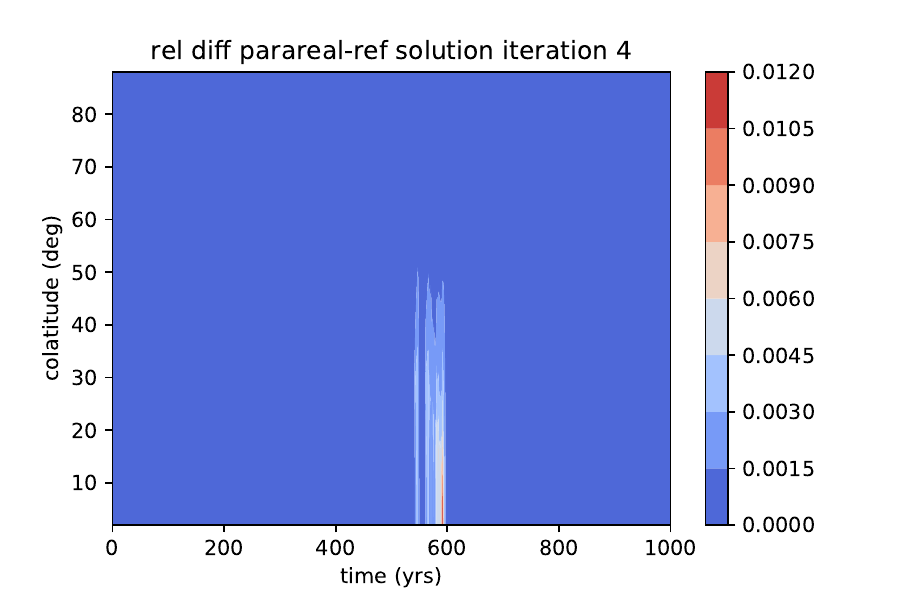}
	\caption{Top: Spatial mean of parareal and reference and macro solution after first (left) and fourth (right) iteration for $N=100$ subintervals, using Euler method in the micro model and constant forcing in the macro model. Bottom: pointwise relative difference between parareal and reference solution.\label{fig:plots}}
\end{figure}

\begin{figure}
	\centering
	\includegraphics[scale=0.475]{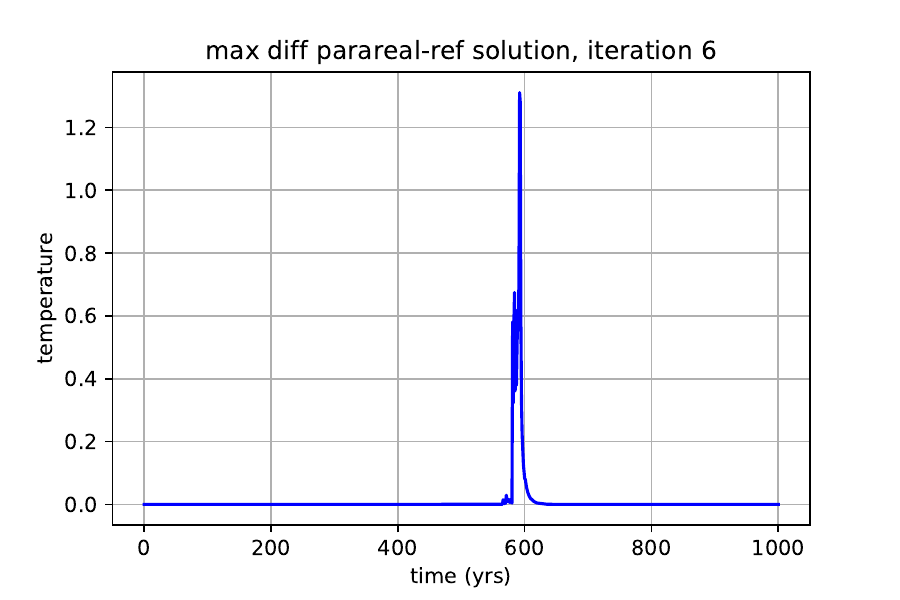}
	\includegraphics[scale=0.475]{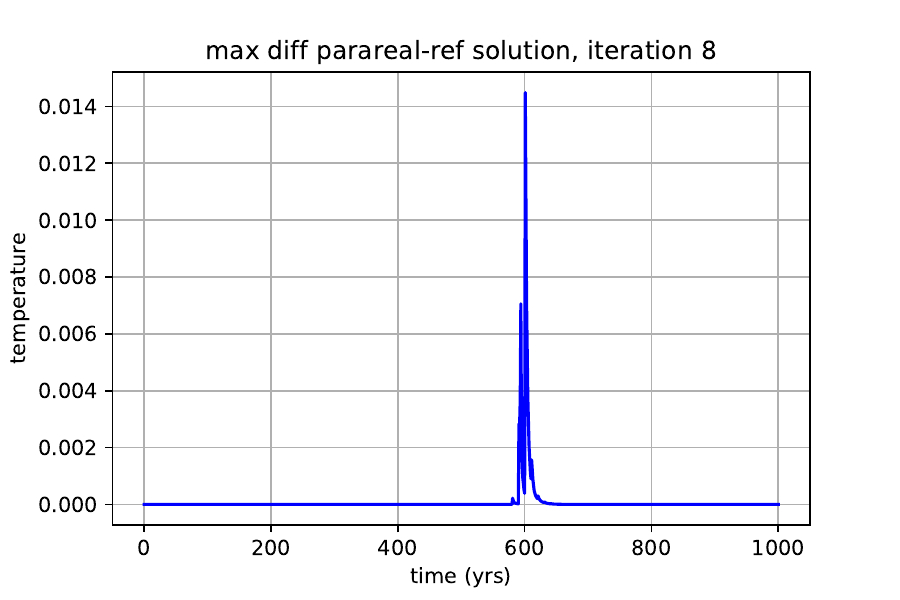}
	\caption{Temporal variation of maximum difference  (in space) between parareal and reference  solution after 6th (left) and 8th (right) iteration for $N=100$ subintervals, using Euler method in the micro model and constant forcing in the macro model. \label{fig:maxerror}}
\end{figure}

Figure \ref{fig:plots} shows as example the solution obtained by the parareal method  for the setting with $N=100$ subintervals (Euler method in the micro, constant forcing in the macro model). This choice of subintervals showed the worst convergence behavior w.r.t. the needed number of parareal iterations, compare Figure \ref{fig:conv-2} and Table \ref{table:1}. It can be seen that even in this case already after 4 iterations there is no  visible difference between parareal and reference solution in the spatial mean, and also only  a small one in the spatial state.
 The main  difficulty for the method is to capture the steep gradient at the time instant close to $t=b_{2}$ (compare \eqref{Qdist2}), which is the second upward jump in solar  forcing.  This can be seen in Figure \ref{fig:maxerror}. It shows the difference between reference and parareal solution over time. Between  6th and 8th iteration, a significant reduction (by a factor $\approx$ 100) can be seen.

No big difference can be seen  in the convergence behavior between the use of both macro model versions. Both versions of the   macro model  give the same precision after nearly the same number of parareal ierations. 
\subsection{Results using adaptive method as micro propagator}
\begin{figure}
	\centering
	\includegraphics[scale=0.3]{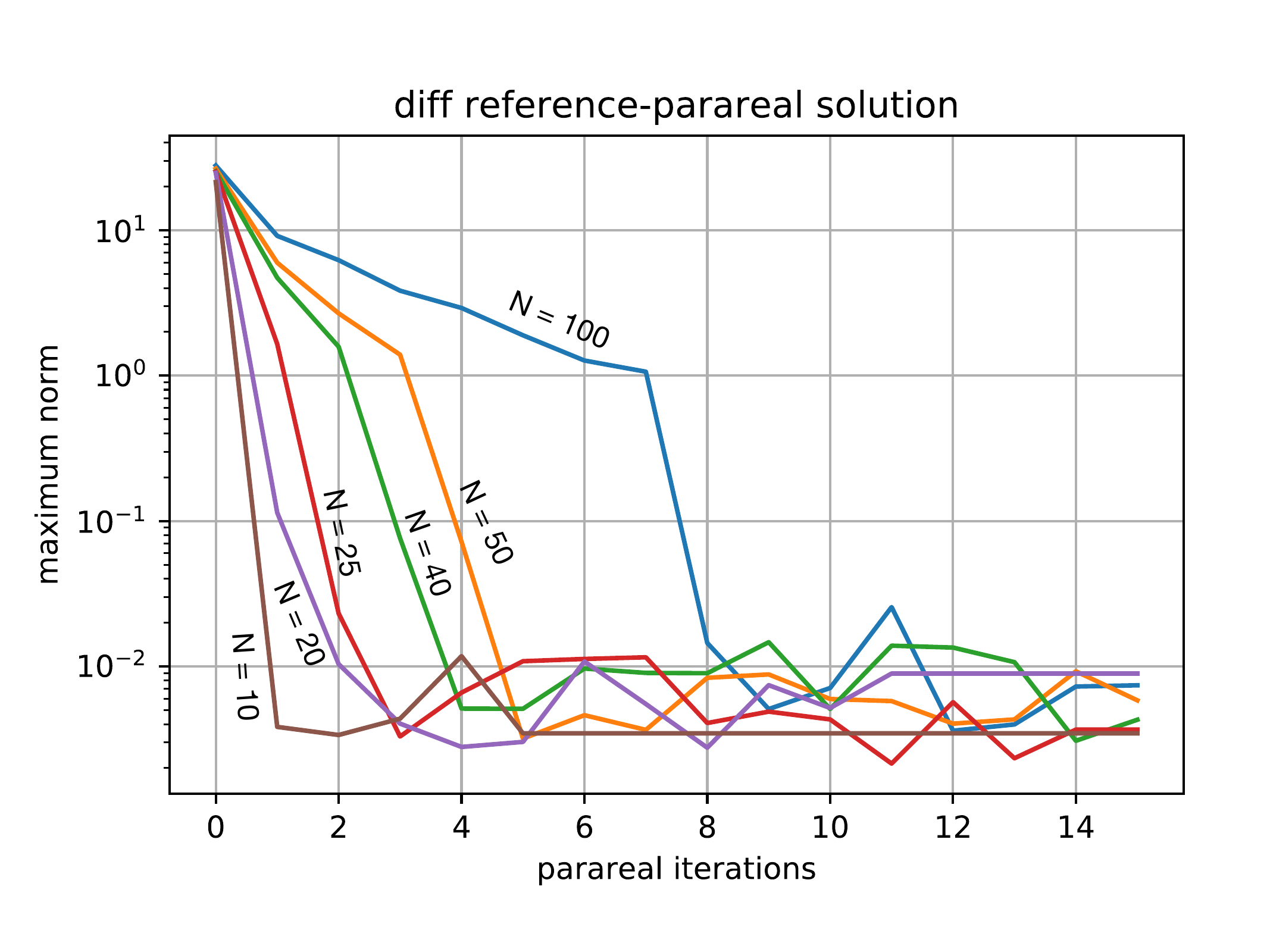}
	\caption{Same as Figure \ref{fig:conv-1}, but  using {\tt lsoda} library routine for the micro model, and the macro model with  varying forcing. \label{fig:conv-3}}
\end{figure}

\begin{figure}
	\centering
	\includegraphics[scale=0.3]{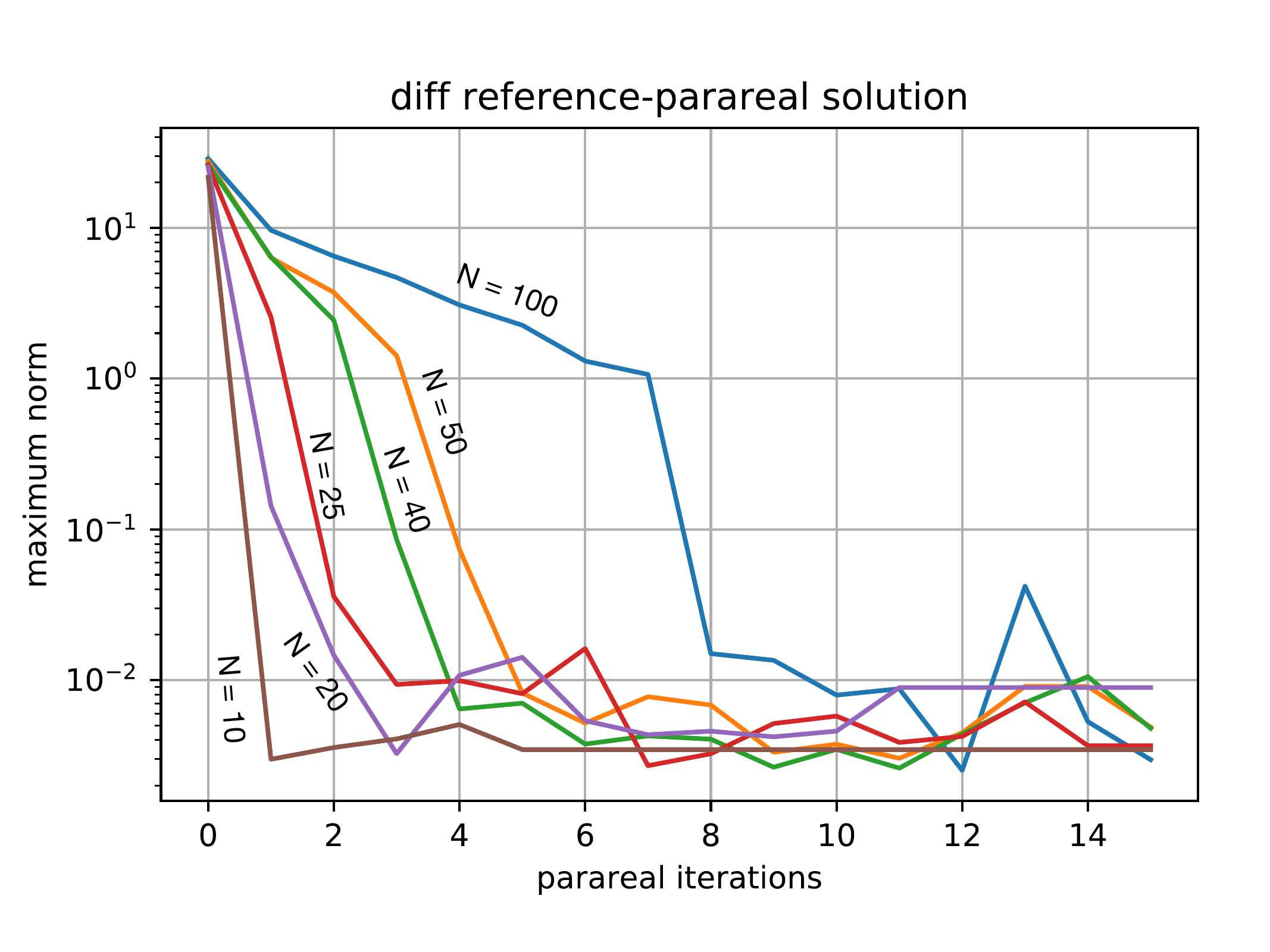}
	\caption{Same as Figure \ref{fig:conv-3}, but   using  the macro model with constant forcing. \label{fig:conv-4}}
\end{figure}

 Figures \ref{fig:conv-3} and \ref{fig:conv-4} show the convergence of the parareal method using the \verb$lsoda$ library routine for the micro model. Here, the choice of the macro model leads to no relevant differences in any case.
Due to the reasons already mentioned at the beginning of this section,  a reduction of the difference between reference and parareal solution up to machine precision could not be reached in this setting.  We have chosen the tolerances of the method such that a difference of less than $10^{-2}$ was reached. This allows us to compare  the results, shown in Table \ref{table:1}, to those of the Euler method. No big difference in the number of iterations needed to give a maximal error less then $10^{-2}$ in both methods can be seen.   


\begin{figure}
	\centering
	\includegraphics[scale=0.465]{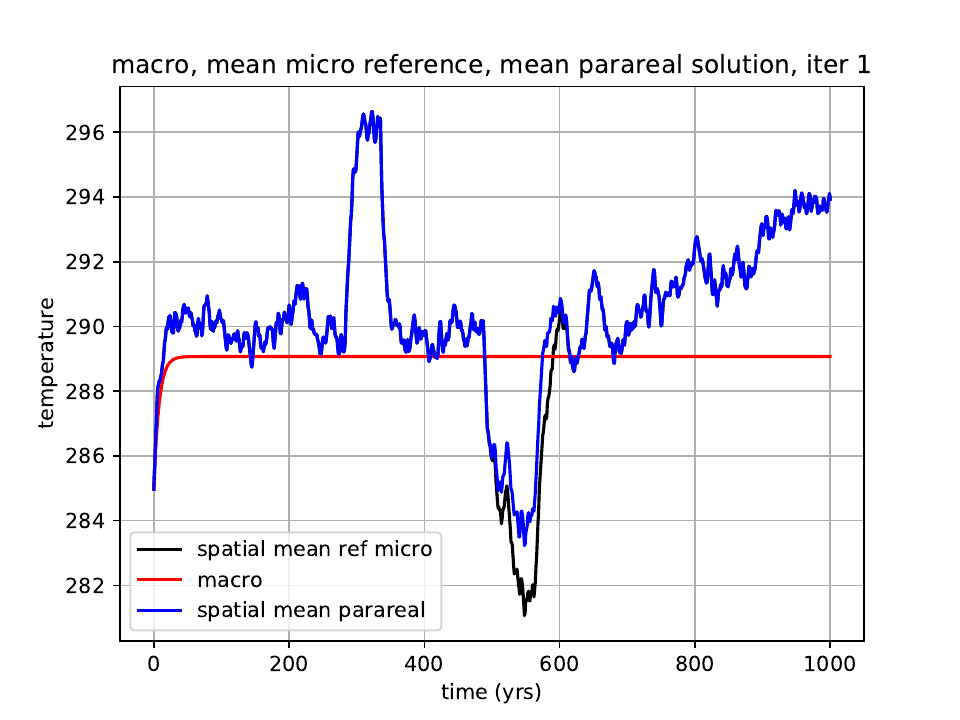}
	\includegraphics[scale=0.465]{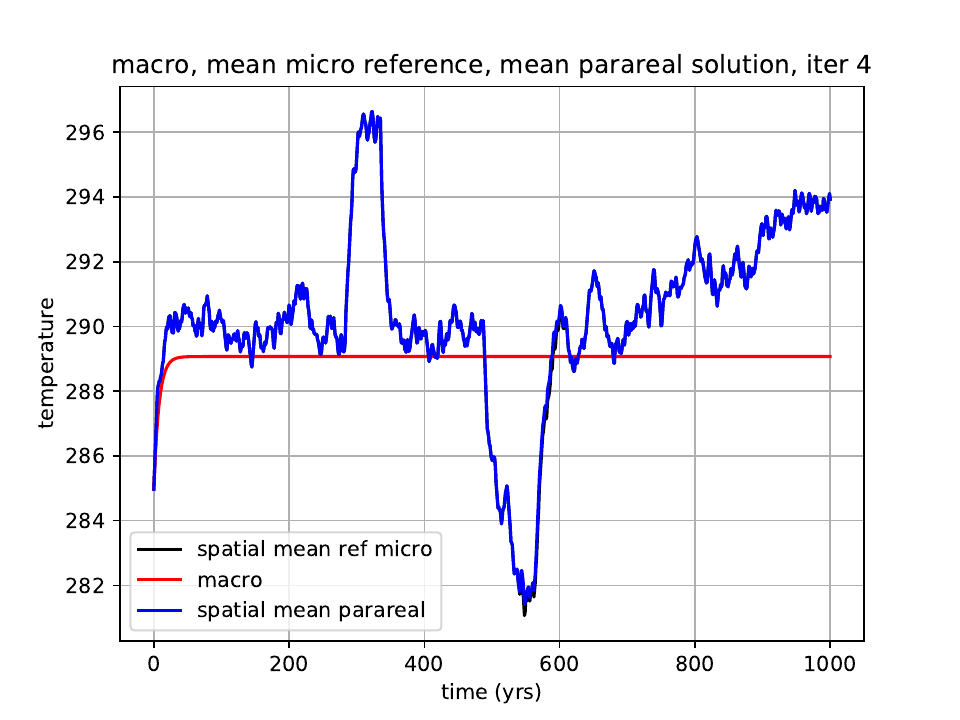}
	\includegraphics[scale=0.475]{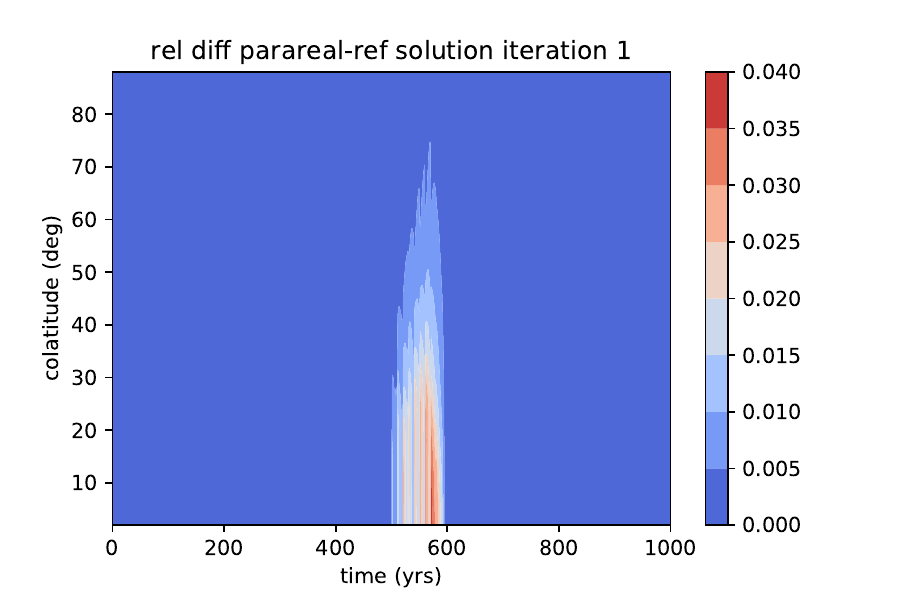}
	\includegraphics[scale=0.475]{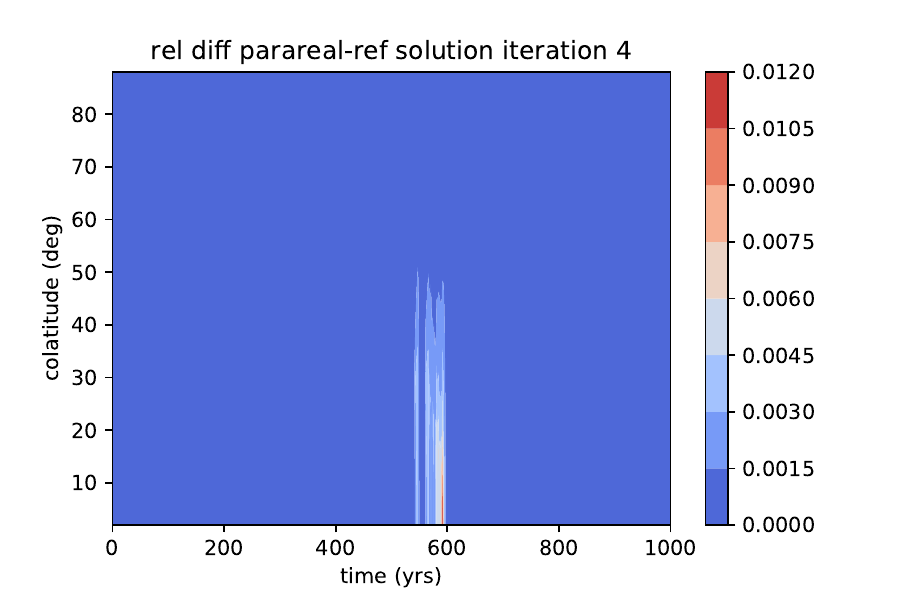}
	\caption{Same as Figure \ref{fig:plots}, but now for the {\tt lsoda} library routine as micro model propagator after the first (left) and the fourth  parareal iteration.\label{fig:plots2}}
\end{figure}

 Figure \ref{fig:plots2} shows again as example the solution obtained  for the setting with $N=100$ subintervals with constant forcing in the macro model. Also here, this choice of subintervals showed the worst convergence  w.r.t. the needed number of parareal iterations, compare Figure \ref{fig:conv-4}. The behavior is comparable to the results obtained using the Euler method: After four iterations no big difference is visible.

\section{Summary and conclusions \label{sec:discussion}}

 We applied the micro/macro parareal method to a 1-D  climate model with  temporally multi-scale forcing. 
As macro model, we use two 0-D versions of the model with spatially averaged coefficients. 

We investigated four configurations (micro model with explicit Euler method with constant stepsize or  adaptive library routine on one hand and macro model with varying or constant forcing on the other) and  a variety of number of parareal subintervals. 

{To estimate the computational gain of the micro-macro method, we considered the number of parareal iterations that were necessary to obtain a reasonable accuracy. This means that we discarded the computational cost of both the macro solver and any communication, and that we assume  an equal distribution of load between subintervals. For the adaptive micro-model, this will not necessarily be the case. As consequence, we  obtained 
only a theoretical bound on computational gain. Concerning the cost of the macro solver, the overhead actually is negligible, since the macro model is a single ODE. The relation between one serial micro and one macro run is about 350:1 (for variable forcing) and up to 1000:1 (for constant forcing) in the used spatial and temporal resolution.}

In all configurations, a reasonable difference to the respective reference solution  (obtained serially with the same time integrator) was obtained after quite a low number of parareal iterations. 
Using the Euler method for the micro propagator, machine precision was reached in {fewer} iterations than numbers of used subintervals. For the adaptive library routine, this was not the case due to the different internal time-grids used  on the whole time interval and when restarted on each subinterval.

There was no relevant difference between the two versions of the macro model: The more inaccurate version with constant forcing was able to predict the solution in all cases but one as well as the one with time-varying forcing.

We conclude that the applied micro/macro parareal algorithm is appropriate for this kind of problem, even if the macro model does not include the temporal multi-scale  features. The {theroretical bound on}  computational gain 
could be even increased when lower (but still reasonable) accuracy requirements are used.
Our results motivate the application of the micro/macro parareal method to more realistic climate models with higher spatial dimension and resolution.

\section*{Acknowledgements}
GS acknowledges the support of the Research Council of the University of Leuven through grant 'PDEOPT', and of the Research Foundation – Flanders (FWO –- Vlaanderen) under grant G.A003.13.
TS acknowledges the support of the German Federal Ministry of Education and Research (BMBF) as part of the Research for Sustainability Initiative (FONA) in the Palmod project,  under grant 01LP1514A.

\bibliographystyle{plain}
\bibliography{parareal_climate}

\begin{thebibliography}{10}

\bibitem{bal2005convergence}
G.~Bal.
\newblock On the convergence and the stability of the parareal algorithm to
  solve partial differential equations.
\newblock In R.~Kornhuber, R.~Hoppe, J.~P{\'e}riaux, O.~Pironneau, O.~Widlund,
  and J.~Xu, editors, {\em Domain decomposition methods in science and
  engineering}, volume~40 of {\em Lecture Notes in Computational Science and
  Engineering}, pages 425--432. Springer Berlin Heidelberg, 2005.

\bibitem{Bay91}
D.~Bayer.
\newblock {\em Einfache mathematische {M}odelle zur {B}eschreibung globaler
  {K}lima{\"a}nderungen}.
\newblock Verlag Dr. Kovac Hamburg, 1991.

\bibitem{BBK}
A.~Blouza, L.~Boudin, and S.-M. Kaber.
\newblock Parallel in time algorithms with reduction methods for solving
  chemical kinetics.
\newblock {\em Communications in Applied Mathematics and Computational
  Science}, 5(2):241--263, 2010.

\bibitem{Brook2008}
E.~Brook.
\newblock Windows on the greenhouse.
\newblock {\em Nature}, 453:291--292, 2008.

\bibitem{engblom2009parallel}
S.~Engblom.
\newblock Parallel in time simulation of multiscale stochastic chemical
  kinetics.
\newblock {\em Multiscale Modeling and Simulation}, 8:46--68, 2009.

\bibitem{Eva98}
L.C. Evans.
\newblock {\em {P}artial {D}ifferential {E}quations}.
\newblock American Math. Society, Providence, Rhode Island, 1998.

\bibitem{farhat2003time}
C.~Farhat and M.~Chandesris.
\newblock Time-decomposed parallel time-integrators: theory and feasibility
  studies for fluid, structure, and fluid--structure applications.
\newblock {\em International Journal for Numerical Methods in Engineering},
  58(9):1397--1434, 2003.

\bibitem{fischer2005parareal}
P.~Fischer, F.~Hecht, and Y.~Maday.
\newblock A parareal in time semi-implicit approximation of the
  {N}avier-{S}tokes equations.
\newblock In R.~Kornhuber, R.~Hoppe, J.~P\'eriaux, O.~Pironneau, O.~Widlund,
  and J.~Xu, editors, {\em Domain decomposition methods in science and
  engineering}, volume~40 of {\em Lecture Notes in Computational Science and
  Engineering}, pages 433--440. Springer Berlin Heidelberg, 2005.

\bibitem{gander2007analysis}
M.J. Gander and S.~Vandewalle.
\newblock Analysis of the parareal time-parallel time-integration method.
\newblock {\em SIAM Journal on Scientific Computing}, 29:556--578, 2007.

\bibitem{Gander2008}
M.~J. Gander and E.~Hairer.
\newblock {Nonlinear {C}onvergence {A}nalysis for the {P}arareal {A}lgorithm}.
\newblock In U.~Langer, M.~Discacciati, D.~E. Keyes, and W.~Zulehner, editors,
  {\em {Domain {D}ecomposition {M}ethods in {S}cience and {E}ngineering
  {XVII}}}, volume~60 of {\em {Lecture Notes in Computational Science and
  Engineering}}, pages 45--56. Springer, 2008.
  
\bibitem{Ganopolski2017}
A.~Ganopolski and V.~Brovkin.
\newblock Simulation of climate, ice sheets and co$_{2}$ evolution during the
  last four glacial cycles with an earth system model of intermediate
  complexity.
\newblock {\em Climate of the Past}, 13(12):1695--1716, 2017.

\bibitem{garrido2005}
I.~Garrido, M.~Espedal, and G.~Fladmark.
\newblock A convergent algorithm for time parallelization applied to reservoir
  simulation.
\newblock In R.~Kornhuber, R.~Hoppe, J.~P{\'e}riaux, O.~Pironneau, O.~Widlund,
  and J.~Xu, editors, {\em Domain Decomposition Methods in Science and
  Engineering}, volume~40 of {\em Lecture Notes in Computational Science and
  Engineering}, pages 469--476. Springer Berlin Heidelberg, 2005.

\bibitem{Ghi75}
M.~Ghil.
\newblock Steady-{S}tate {S}olutions of a {D}iffusive {E}nergy-{B}alance
  {C}limate {M}odel and {T}heir {S}tability.
\newblock Technical Report IMM 410, New York Unversity, Courant Institute of
  Mathematical Sciences, May 1975.

\bibitem{GhiChi87}
M.~Ghil and S.~Childress.
\newblock {\em Topics in {G}eophysical {F}luid {D}ynamics: {A}tmospheric
  {D}ynamics, {D}ynamo {T}heory, and {C}limate {D}ynamics}, volume~60 of {\em
  Applied Math. Sci.}
\newblock Springer, 1987.

\bibitem{Hin83}
A.~C. Hindmarsh.
\newblock {ODEPACK}, {A} {S}ystematized {C}ollection of ode solvers.
\newblock In R.~S.~Stepleman et~al., editor, {\em Scientific Computing},
  volume~1 of {\em IMACS Transactions on Scientific Computation}, pages 55--64.
  North-Holland, 1983.

\bibitem{keller1968numerical}
H.B. Keller.
\newblock {\em Numerical methods for two-point boundary-value problems}.
\newblock Blaisdell (Waltham, MA), 1968.

\bibitem{LeLeSa13}
F.~Legoll, T.~Leli\`{e}vre, and G.~Samaey.
\newblock A micro-macro parareal algorithm: Application to singularly perturbed
  ordinary differential equations.
\newblock {\em S{IAM} {J}. {S}ci. {C}omput.}, 35(4):1951--1986, 2013.

\bibitem{LiMaGa01}
J.-L. {Lions}, Y.~{Maday}, and G.~{Turinici}.
\newblock {R\'esolution d'EDP par un sch\'ema en temps ``parar\'eel''.}
\newblock {\em {C. R. Acad. Sci., Paris, S\'er. I, Math.}}, 332(7):661--668,
  2001.

\bibitem{lubich1987multi}
C.~Lubich and A.~Ostermann.
\newblock Multi-grid dynamic iteration for parabolic equations.
\newblock {\em BIT Numerical Mathematics}, 27(2):216--234, 1987.

\bibitem{maday41parareal}
Y.~Maday.
\newblock Parareal in time algorithm for kinetic systems based on model
  reduction.
\newblock In A.~Bandrauk, M.C. Delfour, and C.~Le~Bris, editors, {\em
  High-dimensional partial differential equations in science and engineering},
  volume~41 of {\em CRM Proceedings and Lecture Notes}, pages 183--194.
  American Mathematical Society, 2007.

\bibitem{maday2002parareal}
Y.~Maday and G.~Turinici.
\newblock A parareal in time procedure for the control of partial differential
  equations.
\newblock {\em Comptes Rendus de l'Acad\'{e}mie des Sciences - Series I -
  Mathematics}, 335(4):387--392, 2002.

\bibitem{maday2005parareal}
Y.~Maday and G.~Turinici.
\newblock The parareal in time iterative solver: a further direction to
  parallel implementation.
\newblock In R.~Kornhuber, R.~Hoppe, J.~P{\'e}riaux, O.~Pironneau, O.~Widlund,
  and J.~Xu, editors, {\em Domain decomposition methods in science and
  engineering}, volume~40 of {\em Lecture Notes in Computational Science and
  Engineering}, pages 441--448. Springer Berlin Heidelberg, 2005.

\bibitem{McgHen14}
K.~McGuffie and A.~Henderson-Sellers.
\newblock {\em The Climate Modelling Primer}.
\newblock Wiley, Chichester, 4th edition, 2014.

\bibitem{mitran2010time}
S.~Mitran.
\newblock Time parallel kinetic-molecular interaction algorithm for {CPU}/{GPU}
  computers.
\newblock {\em Procedia Computer Science}, 1:745--752, 2010.

\bibitem{nievergelt1964parallel}
J.~Nievergelt.
\newblock Parallel methods for integrating ordinary differential equations.
\newblock {\em Communications of the ACM}, 7(12):731--733, 1964.

\bibitem{Nor75}
G.~R. North.
\newblock Theory of {E}nergy-{B}alance {C}lmeate {M}odels.
\newblock {\em J. Atmos. Sci.}, 32:2033--2043, 1975.

\bibitem{Pet83}
L.~R. Petzold.
\newblock Automatic selection of methods for solving stiff and nonstiff systems
  of ordinary differential equations.
\newblock {\em S{IAM} {J}. {S}ci. {C}omput.}, 4:136--148, 1983.

\bibitem{staff2005stability}
G.~Staff and E.~R{\o}nquist.
\newblock Stability of the parareal algorithm.
\newblock In R.~Kornhuber, R.~Hoppe, J.~P{\'e}riaux, O.~Pironneau, O.~Widlund,
  and J.~Xu, editors, {\em Domain decomposition methods in science and
  engineering}, volume~40 of {\em Lecture Notes in Computational Science and
  Engineering}, pages 449--456. Springer, 2005.

\bibitem{Sto14}
T.~Stocker.
\newblock {\em Introduction to {C}limate {M}odeling ({L}ecture {N}otes)}.
\newblock Physics Institute, University of Bern, 2014.

\bibitem{Tro10}
F.~Tr{\"o}ltzsch.
\newblock {\em Optimal Control of Partial Differential Equations: Theory,
  Methods and Applications}.
\newblock Graduate Studies in Mathematics. American Mathematical Society, 2010.

\bibitem{vandewalle1992efficient}
S.~Vandewalle and R.~Piessens.
\newblock Efficient parallel algorithms for solving initial-boundary value and
  time-periodic parabolic partial differential equations.
\newblock {\em SIAM Journal on Scientific and Statistical Computing},
  13:1330--1346, 1992.

\end{thebibliography}

\appendix

\section{Discretization of equation~\eqref{eq:1}\label{sec:discr}}

\subsection{Spatial discretization}
In this appendix, we treat the spatial discretization of the right-hand side of \eqref{eq:1}, which basically means the discretization of the diffusion term  \eqref{eq:diffterm}. Here, we omit the temporal dependency of $T$ in the notation. All indices referring to spatial discretization are denoted as subscripts, whereas indices coming from the temporal discretization will (later on) be denoted as superscripts.

We introduce a spatial grid for $\phi\in[0,\pi/2]$ with stepsize $\Delta\phi=\pi/(2I),I\in\mathbb{N}$.
The gridpoints then are $\phi_{i}=i\Delta\phi,i=0,\ldots,I$. On these points we compute the approximate  solution denoted by 
$\bT_{i}\approx T(\phi_{i})$ for every discrete time-step. We also need intermediate points
$\phi_{i+\frac12}=(i+\frac12)\Delta\phi,i=0,\ldots,I-1$.

We discretize the diffusion term \eqref{eq:diffterm}
by applying central finite differences with stepsize $(\Delta\phi)/2$ twice:
We approximate the inner derivative on the intermediate grid points by
\begin{eqnarray}
\label{eq:theta}
\frac{dT}{d \phi}(\phi_{i+\frac12})&\approx&\frac{\bT_{i+1}-\bT_{i}}{\Delta\phi}
=:\Theta_{i+\frac12},\quad i=0,\ldots,I-1.
\end{eqnarray}
Then, the  derivative  of 
\begin{eqnarray*}
 F(\phi)
&:=&k(\phi,T(\phi))\sin (\phi)\frac{\partial T}{\partial \phi}(\phi)
\end{eqnarray*}
at the points $\phi_{i}$ is  computed,  again by central finite differences: 
\begin{eqnarray}
\label{eq:4}
\left.\nabla\cdot (k{(\phi,T(\phi))}\nabla T(\phi))\right|_{\phi=\phi_{i}}=\frac{dF}{d \phi}(\phi_{i})
\approx
\frac{F(\phi_{i+\frac12})-F(\phi_{i-\frac12})}{\Delta\phi}, 
\quad i=1,\ldots,I-1.
\end{eqnarray}
We need the values of the function $k=k(\phi,T(\phi))$ on the intermediate points:
\begin{eqnarray*}
\kappa_{i+\frac12}&:=&k(\phi_{i+\frac12},T(\phi_{i+\frac12}))
\;\approx \;k\left(\phi_{i+\frac12},\frac{\bT_{i}+\bT_{i+1}}{2}\right),\quad i=0,\ldots,I-1,
\end{eqnarray*}
where the unknown temperature values  at the intermediate points are interpolated linearly.
For the outmost points $\phi_{\frac12},\phi_{I-\frac12}$ this gives, using the homogenous Neumann boundary conditions:
\begin{eqnarray*}
\kappa_{\frac12}
\approx k(\phi_{\frac12},\bT_{1}),&&
\kappa_{I-\frac12}
\approx k(\phi_{I-\frac12},\bT_{I-1}).
\end{eqnarray*}
Eventually, needed values of the coefficient $k$ (if not constant) at the intermediate points $\phi_{i+\frac12}$  have to be interpolated accordingly. 
We then obtain for the  terms on the right-hand side of  \eqref{eq:4}:
\begin{eqnarray}
\label{eq:5}
F(\phi_{i+\frac12})
\approx
\kappa_{i+\frac12}\sin (\phi_{i+\frac12})\Theta_{i+\frac12},&&
F(\phi_{i-\frac12})\approx\kappa_{i-\frac12}\sin (\phi_{i-\frac12})\Theta_{i-\frac12}.
\end{eqnarray}

\subsection{Time discretization}
We now obtain a nonlinear nonautonomous system of ODEs
\begin{eqnarray}
\bT'(t)&=&f(t,\bT(t)),\quad t\ge0\label{eq:micro_app}
\end{eqnarray}
with $\bT(t),f(t,\bT(t))\in \mathbb{R}^{I-1}$ and some initial value $\bT(0)\in \mathbb{R}^{I-1}$. The $i$-th component of the function $f$ is given by \eqref{eq:1} using \eqref{eq:theta}, \eqref{eq:4}, and \eqref{eq:5}.
For both models, we re-scaled the time from seconds  to years, i.e. we set
\begin{eqnarray*}
\tilde t&:=& \frac{t}{s_{year}}\quad \text{with } s_{year}=60\times60\times24\times 365 = 3.1536\times10^{7}.
\end{eqnarray*}
The time derivative on the left-hand side of both models can be transformed using the definition $\tilde T(\tilde t):=T(s_{year}\tilde t)$ and the formula
\begin{align*}
\frac{dT}{dt}(t)=\frac{dT}{d(s_{year}\tilde t)}(s_{year}\tilde t)=
\frac{1}{s_{year}}\frac{dT}{d\tilde t}(s_{year}\tilde t)
=\frac{1}{s_{year}}\frac{d\tilde T}{d\tilde t}(\tilde t).
\end{align*}
The dependency on $\phi$ in the 1-D model was suppressed here for simplicity.
Then  both models, \eqref{eq:0-D} and \eqref{eq:1}, retain their formulation (omitting the tildes) when 
the respective right-hand sides are multiplied by $s_{year}$.

\end{document}